\documentclass[12pt]{article}
\usepackage{color}
\usepackage{amsmath, amsthm, amsfonts, amssymb,latexsym}
\usepackage{graphicx} 
\usepackage{epstopdf}
\usepackage[bookmarks=false]{hyperref}
\usepackage{mathrsfs}
\UseRawInputEncoding

\newtheorem{thm}{Theorem}[section]
\newtheorem{cor}[thm]{Corollary}
\newtheorem{lem}[thm]{Lemma}

\theoremstyle{definition}
\newtheorem{defn}[thm]{Definition}

\newtheorem{exam}[thm]{Example}

\begin{document}

\def\prof{{\sc Proof.\ \ }}
\def\sect#1{\begin{center}\section{#1}\end{center}} 
\def\R#1{{\bf R}^{#1}}
\def\I#1#2{\int\limits_{#1}#2} 
\def\p#1{#1^{\prime}} 
\def\l#1#2#3{\lim_{#1 \rightarrow #2}#3}
\def\ld#1#2#3{\liminf_{#1 \rightarrow #2}#3} 
\def\lu#1#2#3{\limsup_{#1 \rightarrow #2}#3}
\def\E#1#2{{\bf E}^{#1}\left(#2\right)} 
\def\P#1#2{{\bf P}^{#1}\left(#2\right)} 
\def\Pt#1#2#3{{\bf P}(#1,#2,#3)} 
\def\pb#1{{\bf P}\left(#1\right)}
\def\eb#1{{\bf E}\left(#1\right)}
\def\Di#1{\,Dim\,#1}
\def\di#1{\,dim\,#1}
\def\fr#1#2{\frac{#1}{#2}}
\def\beq{\begin{equation}}
\def\eeq{\end{equation}}
\def\bea{\begin{eqnarray}}
\def\bean{\begin{eqnarray*}}
\def\eean{\end{eqnarray*}}
\def\eea{\end{eqnarray}}
\def\heq#1#2#3{\hbox to \hsize{\hskip #1 $#2$ \hss (#3)}}
\def\df#1#2{\frac{\displaystyle #1}{\displaystyle #2}}

\def\bdes{\begin{description}}
\def\ndes{\end{description}}

\newcommand{\bh}{{\bf h}}
\newcommand{\hf}{{\bf f}}
\newcommand{\he}{{\bf e}}
\newcommand{\hL}{{\bf L}}
\newcommand{\hg}{{\bf g}}
\newcommand{\hG}{{\bf G}}
\newcommand{\hM}{{\bf M}}

\newcommand{\bb}{\mathbb{b}}
\newcommand{\ww}{\mathbb{W}}
\newcommand{\hh}{\mathbb{H}}
\newcommand{\dd}{\mathbb{D}}
\newcommand{\cc}{\mathbb{C}}
\newcommand{\ee}{\mathbb{E}}
\newcommand{\zz}{\mathbb{Z}}
\newcommand{\nn}{\mathbb{N}}
\newcommand{\pp}{\mathbb{P}}
\newcommand{\qq}{\mathbb{Q}}
\newcommand{\ttt}{\mathbb{T}}

\def\rr{\Bbb R}
\def\a{\alpha}
\def\b{\beta}
\def\g{\gamma}
\def\s{\sigma}
\def\ep{\epsilon}
\def\d{\delta}
\def\D{\Delta}
\def\o{\omega}
\def\O{\Omega}
\def\n={\not=}
\def\u>{\wedge}
\def\d>{\vee}
\def\.{\bullet}
\def\l{\lambda}
\def\L{\Lambda}
\def\r{\rho}
\def\vf{\varphi}
\def\f{\phi}
\def\t{\tau}
\def\z{\zeta}
\def\na{\nabla}

\def\B{{I\!\!B}}
\def\E{{I\!\!E}}
\def\N{{I\!\!N}}
\def\P{{I\!\!P}}
\def\Q{{I\!\!\!Q}}
\def\R{{I\!\!R}}

\def\AA{\mathcal A}
\def\BB{\mathcal B}
\def\CC{\mathcal C}
\def\DD{\mathcal D}
\def\FF{\mathcal F}
\def\EE{\mathcal E}
\def\GG{\mathcal G}
\def\JJ{\mathcal J}
\def\LL{\mathcal L}
\def\NN{\mathcal N}
\def\PP{\mathcal P}
\def\SS{\mathcal S}

\def\Fht{\hbox{${\cal F}_{t}$}}
\def\Ght{\hbox{${\cal G}_{t}$}}

\def\U{\bigcup}
\def\Uu{\bigcap}
\def\Au{\forall}
\def\Eu{\exists}
\def\8u{\infty}
\def\0{\emptyset}
\def\mp{\longmapsto}
\def\rt{\rightarrow}
\def\up{\uparrow}
\def\dn{\downarrow}
\def\sub{\subset}
\def\vep{\varepsilon}

\def\ln{\langle}
\def\rn{\rangle}

\def\<<{\langle\!\langle}
\def\>>{\rangle\!\rangle}
\def\3|{|\!|\!|}

\def\ess{\text{\rm{ess}}}
\def\beg{\begin}

\def\beqt{\begin{equation}}
\def\neqt{\end{equation}}
\def\beq{\begin{equation}}
\def\neq{\end{equation}}

\def\bdes{\begin{description}}
\def\ndes{\end{description}}

\def\Ric{\text{\rm{Ric}}}
\def\Hess{\text{\rm{Hess}}}
\def\i{\text{\rm{i}}}
\def\ii{\text{\rm{ii}}}
\def\iii{\text{\rm{iii}}}
\def\iv{\text{\rm{iv}}}
\def\v{\text{\rm{v}}}
\def\vi{\text{\rm{vi}}}
\def\vii{\text{\rm{vii}}}
\def\viii{\text{\rm{viii}}}
\def\e{\text{\rm{e}}}

\def\Ra{\Rightarrow}
\def\Lra{\Leftrightarrow}
\def\rto{\longrightarrow}
\def\la{\leftarrow}
\def\ra{\rightarrow}
\def\ua{\uparrow}
\def\da{\downarrow}


\title{Multiparametric analysis of conic linear optimization based on lift-and-project procedure
\thanks{Supported by the National Natural Science Foundation of China (11871118,12271061).}}

\author{Zi-zong Yan \thanks{School of Information and Mathematics,
Yangtze University, Jingzhou, Hubei,
China(zzyan@yangtzeu.edu.cn).}
\ Xiang-jun Li \thanks{School of Information and Mathematics,
Yangtze University, Jingzhou, Hubei,
China(franklxj001@163.com).}
\ and Jin-hai Guo\thanks{School of Information and Mathematics,
Yangtze University, Jingzhou, Hubei,
China(xin3fei@21cn.com).}}
\date{}
\maketitle

\begin{abstract} We study the application of the lift-and-project procedure to the multiparametric analysis of conic linear optimization (CLO) problems. We first introduce the concept of a pair of primal and dual conic representable sets, and define the set-valued mappings between them. We then explore a novel duality of mpCLOs, that allows us to generalize and treat the previous results for parametric analysis in a unified framework. In particular, we discuss the behavior of the optimal partition of a conic representable set. This leads to invariant region decomposition of a conic representable set, which is more general than the results in the literature. Finally, we study the properties of the optimal objective values as a function of parametric vectors. All the results are corroborated by examples with correlations.

\textbf{Keywords:} Multiparametric conic linear optimization, conic representable set, lift-and-project, optimal partition, duality, multiparametric KKT property \\

\textbf{AMS subject classifications.} Primary: 90C31; Secondary: 90C25, 90C05, 90C22, 90C46
\end{abstract}

\section{Introduction}
We are interested in the following multiparametric conic linear optimization (mpCLO) problems with two independent vectors of parameters $u,v \in \mathbb{R}^r$
\begin{equation} \label{primaljihe1} \tag{P-lift} \begin{array}{lll} p^*(u)=& \min & \langle c+M^Tu,x \rangle \\ & s.t. & Ax=b,\\ & & x \in K \end{array} \end{equation}
and
\begin{equation} \label{primaljihe2} \tag{D-lift} \begin{array}{lll} d^*(v)=& \min & \langle d+M^Tv,y \rangle \\ & s.t. & By=a,\\ && y \in K^* \end{array} \end{equation} for given vectors $a \in \mathbb{R}^l$, $b \in \mathbb{R}^m$, $c,d \in \mathbb{R}^q$ and matrices $A \in \mathbb{R}^{m\times q}$, $B \in \mathbb{R}^{l\times q}$, $M\in \mathbb{R}^{r\times q}$, where $K\subset \mathbb{R}^q$ is a pointed, closed, convex, solid (with a non-empty interior) cone (for this formulation of a primal-dual pair and its properties, Ref. \cite{NN94}); and $K^*$ is the dual of $K$ under the standard inner-product, that is, \[ K^*=\{y\in \mathbb{R}^q| \langle y,x\rangle\geq 0, \forall x\in K\}. \] The minimum values of the objective functions are denoted by $ p^*(u)$ and $d^*(v)$, respectively. These two problems are multiparametric linear programming (mpLP) if $K$ is a nonnegative orthant; and multiparametric semidefinite programming (mpSDP) if $K$ is the cone of symmetric positive semidefinite matrices (for example, Refs. \cite{AL12,BN01,BPT13,Nem07}).

Two technical claims are assumed throughout this paper.

{\bf Assumption 1}. The three range spaces $R(A^T),R(B^T),R(M^T)$ are orthogonal to each other, and their direct sum is equal to the entire space $R^q$.

{\bf Assumption 2}. $MM^T=I_r$, a $r\times r$ unit matrix.

Assumption 1 is necessary for (\ref{primaljihe1}) and (\ref{primaljihe2}) because it guarantees that the perturbation of the objective function is independent of the constraints; More importantly, it guarantees that the (nonstandard) Lagrange dual of one is a projection of the other. We comment more on this in Subsections \ref{nons1} and \ref{otjg1}.

\subsection{Motivations}
The original motivation for this study is the inductive proof of the strong duality theorem of CLOs (see Theorem \ref{dualityt1} and its proof at the end of Subsection \ref{otjg1}). Indeed, such a proof requires that a pair of primal-dual optimal solutions remains the same in the lift-and-project procedure, in which the complement slackness property always holds. More precisely, a pair of optimal solutions $(x^*(u),y^*(v))$ of (\ref{primaljihe1}) and (\ref{primaljihe2}) satisfies the following equation \begin{equation} \label{slack1} \langle x^*(u),y^*(v)\rangle = 0, \end{equation} (see Theorem \ref{interior2yy}). A surprising aspect of this result is that it holds for some vector pairs $(u,v) \in \mathbb{R}^r\times \mathbb{R}^r$, under Slater's type condition.

We exploit the special structure of the CLO problem to obtain a new algebraic representation of the Lagrangian dual problem (see Subsection \ref{nons1}). Such a dual, called a nonstandard dual, results in a unified representation of the feasible sets. Surprisingly, for (\ref{primaljihe1}) and (\ref{primaljihe2}), the nonstandard dual of one is a projection of the other; and their projections also satisfy the complement slackness property \begin{equation} \label{slack2} \langle \bar{x}^*(v),\bar{y}^*(u)\rangle = 0 \end{equation} under the same Slater's type condition, where $(\bar{x}^*(v),\bar{y}^*(u))$ denotes a pair of optimal solutions of the projections (see Theorem \ref{interior2yy}).

These characterization results directly motivate us to explore in detail the relationship between the following two conic (linear inequality) representable sets, \begin{subequations}
\begin{align} \label{primalcone1} \varTheta_P &= \{v \in\mathbb{R}^r| d+M^Tv+B^T w^1\in K \ for\ some\ w^1 \in\mathbb{R}^l \}, \\ \label{primalcone2} \varTheta_D&= \{u \in\mathbb{R}^r| c+M^Tu+A^Tw^2 \in K^* \ for \ some\ w^2 \in\mathbb{R}^m \}, \end{align} \end{subequations} (see Subsection \ref{jzyss}). This relationship helps us present a geometric framework that unifies and extends some properties of parametric LPs and SDPs to the case of mpCLOs.

\subsection{Related works}
Prior to parametric analysis in recent years, the actual invariancy region played an important role in the development of parametric LPs and SDPs. Adler and Monteiro \cite{AD92} investigated the parametric analysis of LPs using the optimal partition approach, in which they identified the range of a single-parameter where the optimal partition remained invariant. Other parametric analysis treatments for LPs based on the optimal partition approach were reported in previous studies \cite{Ber97,Deh07,GG08,Gre94,Ha10,RT05,JR03}. The actual invariancy region has been studied extensively in the SDP setting (Refs. \cite{GS99,MT20}), the second-order conic optimization (Refs. \cite{MT21}), and more generally in CLO (Refs. \cite{Yil04}).

Gass and Saaty \cite{GS55} proposed the first method for solving parametric LPs, and Gal and Nedoma \cite{GN72} generalized this method. Recently, there has been growing interest in multiparametric optimizations arising from process engineering, such as process design, optimization, and control. The survey by Pistikopoulos et al.\cite{Pis12} contributes to recent theoretical and algorithmic advances, and applications in the areas of multi-parametric programming, specifically explicit/multi-parametric model predictive control. Thus far, various types of invariances in parametric/multiparametric problems have been used, mainly in single-parametric or bi-parametric analyses. To the best of our knowledge, almost all approaches to multiparametric LPs reported in the literature (Refs \cite{Bor03,Fil97,Fil04,Gal95,GG97,Sch87,WW90}) use bases to obtain a description of the invariancy regions.

\subsection{Contributions}
This study investigates multiparametric optimization in general CLO problems in which either the objective function or the right-hand side is perturbed along many fixed directions. First, we establish the connection by showing that the conic representable set defined by (\ref{primalcone1}) or (\ref{primalcone2}) is viewed as a projection of the feasible set of an appropriately defined mpCLO problem, although this connection is not easily identifiable. We then develop classical duality theory using set-valued mappings to relate the two conic representable sets. This treatment makes it possible to combine some known, yet scattered, results and derive new ones. All these results can be used to develop the optimal partition approach given in \cite{AD92} for parametric LPs and \cite{GG08} for parametric SDPs.

The main contributions of this study are as follows:

(1) Characterization of the relationship between the primal and dual conic representable sets.

(2) Presentation of a novel duality in CLO.

(3) Identification of optimal partitions and development of parametric analysis technique.

In the first category, we define set-valued mappings between the primal and dual conic representable sets, which provide an important and useful tool. This tool plays a critical role in the analysis.

In the second category, we develop the classical duality theory of CLO. In addition to the complementary slackness properties mentioned earlier, we present the weak and strong duality properties for the pair of almost primal-dual problems (\ref{primaljihe1}) and (\ref{primaljihe2}). We show that the sum of their objective optimal values is a bilinear function with respect to the vectors of parameters $u$ and $v$ (see Corollary \ref{sumtwop21}). The corresponding multiparametric KKT (mpKKT) conditions are also given in Theorem \ref{interior2yy}.

In the third category, along with invariant region decomposition, our results capture and generalize the SDP cases of Mohammad-Nezhad and Terlaky \cite{MT20}. We develop the concepts of the nonlinearity region and the transition point for the optimal partition to conic representable sets by means of set-valued mappings, and provide some sufficient conditions for the existence of a nonlinearity region and a transition point. Such concepts are very useful for the analysis of an mpCLO problem because the nonlinearity region can be regarded as a stability region and its identification has a significant influence on the post-optimal analysis of SDPs.

In addition, we study the behavior of the optimal value function in its domain; and partially answer the open question proposed by Hauenstein et al. \cite{Hau19}.

\subsection{Organization of the paper}
In the next section, we review some useful results from convex analysis and the duality theory of CLO, and introduce nonstandard Lagrange dual and lift-and-project procedures. In Section 3, we define the set-valued mappings between the aforementioned conic representable sets and develop duality theory in CLOs. In Section 4, some fundamental concepts are introduced and several examples are presented to illustrate these concepts. The extension of the corresponding approach to the multiparametric analysis of CLOs is discussed in Section 5. The identification of the optimal partitions of a conic representable set is discussed and the behavior of the optimal value function in its domain is studied in this section. Finally, in Section 6, we conclude the paper with some remarks.

\section{Preliminaries} \label{sect21}
The aim of this section is twofold. First, we introduce the notation used in this study. However, we state some useful tools and results that facilitate the subsequent proofs.

\subsection{Notation}
First, we review some useful facts regarding convex sets and cones. A standard reference for convex analysis is a book by Rockafellar \cite{Roc70}.

The (topological) boundary of a set $C\in\mathbb{R}^q$ is denoted as $\partial C$ and defined as $\partial C =\operatorname{cl}(C) \backslash \operatorname{int}(C)$, where cl($C$) and int($C$) denote the closure and interior of $C$, respectively. $\operatorname{dim}(C)$ denotes the affine dimension of $C$. A set $C$ is called simply connected if for any two points $x,y\in C$, there is a continuous curve $\Gamma\subset C$ connecting $x$ and $y$. A singleton set is referred to as simply connected set. A set $C$ is called convex if for any $x,y\in C$, the linear segment $[x,y]=\{\alpha x+(1-\alpha)y|\alpha\in[0,1]\}$ is contained in $C$. The convex hull of $C$, denoted as $\operatorname{conv}(C)$, is a set of all the convex combinations from $C$.

Let $C$ be a nonempty convex set in $\mathbb{R}^q$. The vector $0\ne h\in \mathbb{R}^q$ is called the recession direction of $C$ if $c+\lambda h \in C$ for all $c\in C$ and $\lambda >0$. The set of all recession directions of $C$ is called the recession cone of $C$, and is denoted as $0^+(C)$. Given a boundary point $\bar{x}$ of $C$, its normal cone $\operatorname{Normal}(C,\bar{x})$ is defined as
\[ \operatorname{Normal}(C,\bar{x})=\{c\in \mathbb{R}^q | \langle c,\bar{x}\rangle\geq \langle c,x\rangle \ \operatorname{for \ all} \ x\in C\}.\]
A point $\bar{x}$ is a vertex of $C$ if its normal cone is full-dimensional.

Mapping $\Phi(\xi): \mathbb{R}^{r} \rightrightarrows \mathbb{R}^{r}$ is called set-valued mapping if it assigns a subset of $ \mathbb{R}^{r}$ to each element of $ \mathbb{R}^{r}$. The following set
\[\Phi(C)=\{\Phi(\xi)|\xi\in C\}\] denotes the image of set $C\in \mathbb{R}^{r}$ under set-valued mapping $\Phi$.

The notations $R(A)=\{Ax|x\in\mathbb{R}^q\}$ and $N(A)=\{x| Ax=0,x\in\mathbb{R}^q\}$ denote the range and kernel space of the matrix $A$, respectively.

Let $U$ be an open subset of $\mathbb{R}^r$. We say that the mapping $f:U\rightarrow \mathbb{R}$ is G$\hat{a}$teaux differentiable at a point $u\in U$ in a direction $h\in \mathbb{R}^r$ if there exists $f'(u,h)\in \mathbb{R}$ such that
\[ f'(u,h)=\lim\limits_{t\rightarrow 0+}\frac{1}{t}(f(u+th)-f(u)),\] where $f'(u, h)$ is called the directional derivative of $f$ at $u$ in the direction $h$.
The mapping $f$ is called G$\hat{a}$teaux differentiable at $u$, provided that $f'(u,h)$ exists for all $h \in \mathbb{R}^r$ and the mapping $f'(u,\cdot):\mathbb{R}^r\rightarrow \mathbb{R}$ is a continuous linear operator.

\subsection{The nonstandard dual} \label{nons1}
Recall that (\ref{primaljihe1}), without perturbation is the typical form of the CLO
\begin{equation} \label{primaljihe10} \tag{P} \begin{array}{ll} \min & \langle c,x \rangle \\ s.t. & Ax=b,\\ & x \in K, \end{array} \end{equation}
and its Lagrangian dual problem is expressed as follows:
\begin{equation} \label{primaljihe20} \tag{D} \begin{array}{ll} \max & \sum\limits_{i=1}^m b_i w_i \\ s.t. & A^Tw \leq_{K^*} c,\\ & w \in\mathbb{R}^m. \end{array} \end{equation}

For clarity and elegance (another crucial reason is given in the next section), we use an equivalent form instead of (\ref{primaljihe20}). In other words, (\ref{primaljihe20}) is rewritten as a new objective function with new constraints. This can be implemented through the following process.

Left multiplying both sides of the equality constraint in (\ref{primaljihe20}) by the matrix $B$ and $M$, respectively, one has $By=Bc$ and $My=Mc$ if $y=c-A^Tw$. Conversely, if $y\in K^*$ satisfies both $By=Bc$ and $My=Mc$, $w\in \mathbb{R}^m$ exists, such as $y=c-A^Tw$. Thus, the feasible set of (\ref{primaljihe20}) can be expressed equivalently as in $\{y\in K^*| By=a, My=Mc\}$. On the other hand, if a pair of vectors $(x;w)$ is feasible, and if $b=Ad$ and $y=c-A^Tw$, \begin{equation}\label{propxz1} \langle d, c-y\rangle =\langle d, A^Tw \rangle = \langle Ad, w\rangle = b^Tw. \end{equation} That is, (\ref{primaljihe20}) can be rewritten as follows:
\begin{equation} \label{primaljihe30} \tag{D-non} \begin{array}{ll} \max & \langle d,c-y \rangle \\ s.t. & By=a,\\ &My=Mc,\\ & y \in K^*, \end{array} \end{equation} where $a=Bc$. This is essentially the dual of (\ref{primaljihe10}), where $y$ denotes the slackness vector variable in (\ref{primaljihe20}). Therefore, (\ref{primaljihe20}) and (\ref{primaljihe30}) are {\it the standard dual and nonstandard dual} of (\ref{primaljihe10}), respectively.

For the primal and dual programs (\ref{primaljihe10}) and (\ref{primaljihe20}), the weak duality property holds. Then, from the equality (\ref{propxz1}), the following result holds:

\cor\label{weakdual1} If $b=Ad$ and $a=Bc$, then for any primal feasible solution $x$ of (\ref{primaljihe10}) and any dual feasible solution $y$ of (\ref{primaljihe30}), the weak duality property holds, i.e.,
\begin{equation} \label{strongdual1} \langle c,x\rangle \geq \langle d,c-y \rangle. \end{equation} The equality holds if and only if $(x,y)$ is a pair of optimal solutions. \upshape

If the primal and dual programs have optimal solutions and the duality gap is zero, that is, the equality (\ref{strongdual1}) holds, then the Karush-Kuhn-Tucker (KKT) conditions for the primal-dual CLO pair are
\begin{subequations}
\begin{align} \label{jbkkt1} & Ax = b,\ x\in K, \\ \label{jbkkt2} & By=a,\ M y=Mc,\ y \in K^*, \\ \label{jbkkt3} & \langle x,y \rangle=0.
\end{align}
\end{subequations} Conversely, if $(x^*,y^*)\in\mathbb{R}^q\times \mathbb{R}^q$ is a pair of solutions of system (\ref{jbkkt1})-(\ref{jbkkt3}), then $(x^*,y^*)$ is a pair of optimal solutions of (\ref{primaljihe10}) and (\ref{primaljihe30}).

\cor\label{weakdual2} The optimal solutions of (\ref{primaljihe10}) and (\ref{primaljihe30}) are independent of the choice of $c$ and $d$ only if the pair $(d,c)$ satisfies $Ad=b$ and $Bc=a$. \upshape

\prof It is a direct consequence of Corollary \ref{weakdual1}. $\square$

In the following statements, we assume $b=Ad$ and $a=Bc$. Alternatively, $d$ and $c$ are feasible for (\ref{primaljihe1}) and (\ref{primaljihe2}), respectively.

We say that a CLO problem is strictly feasible or satisfies Slater's condition if a feasible interior exists. For example, for (\ref{primaljihe10}), a feasible solution, $x$, is strictly feasible if $x\in \operatorname{int}(K)$; whereas for (\ref{primaljihe30}), a feasible solution, $y$, is strictly feasible if $y\in \operatorname{int}(K^*)$. The following strong duality theorem is fundamental (Refs. \cite{Bar02,BS00,BL05,BV04,LV16,Tun11}).

\thm\label{dualityt1} Consider the primal-dual CLO pair (\ref{primaljihe10})-(\ref{primaljihe30}).

(1) If the dual problem is bounded from above and if it is strictly feasible, then the primal problem attains its minimum, and there is no duality gap.

(2) If the primal problem is bounded from above and if it is strictly feasible, then the dual problem attains its maximum, and there is no duality gap. \upshape

\subsection{The lift-and-project procedure}\label{lift1}
Under Assumptions 1 and 2, as discussed in the previous subsection, the nonstandard dual of (\ref{primaljihe1})
can be defined by
\begin{equation} \label{primaljihe3} \tag{P-pro} \begin{array}{lll} \bar{d}^*(u)=& \max & \langle d,c+M^Tu-y \rangle \\ & s.t. & By=a,\\ & & My=Mc+u, \\ && y \in K^*, \end{array} \end{equation} which can be viewed as a projected model of (\ref{primaljihe2}) on the affine set $\{y\in\mathbb{R}^q|My=Mc+u\}$. However, a projected model of (\ref{primaljihe1}) on affine set $\{x\in\mathbb{R}^q|Mx=Md+v\}$ can be described as
\begin{equation} \label{primaljihe4} \tag{D-pro} \begin{array}{lll}\bar{p}^*(v)=& \max & \langle c, d+M^Tv-x \rangle \\ & s.t. & Ax=b,\\ && Mx=Md+v, \\ && x \in K, \end{array} \end{equation} which is the nonstandard dual of (\ref{primaljihe2}). Here the maximum values of the objective functions are denoted as $\bar{d}^*(u)$ and $\bar{p}^*(v)$, respectively.

\defn\label{amlost1} We say that a pair of mpCLOs is almost dual if the nonstandard Lagrange dual of the one is either a projection or a lifting of the other. \upshape

\begin{figure}[htbp]
\centering\includegraphics[width=4.5in]{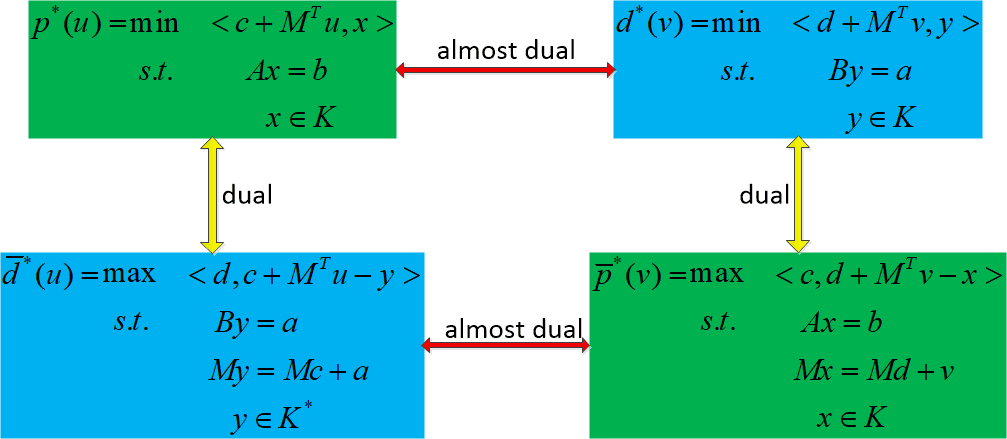}
\caption{Dual and almost dual}\label{fig:1}
\end{figure}

By definition \ref{amlost1}, either (\ref{primaljihe1}) and (\ref{primaljihe2}) or (\ref{primaljihe3}) and (\ref{primaljihe4}) are almost dual. If the objective functions of (\ref{primaljihe3}) and (\ref{primaljihe4}) are replaced by $\langle d,c-y \rangle$ and $\langle c,d-x \rangle$, respectively, then their optimal solutions remain unchanged. Thus, the perturbations in (\ref{primaljihe3}) and (\ref{primaljihe4}) occur on the right-hand side and are not in the objective function data. The main duality examined in this study is shown in Figure \ref{fig:1}.

It is necessary to use the nonstandard dual model because a uniform algebraic representation of the feasible sets helps us to understand the lift-and-project procedure correctly.

\cor \label{xinzd1} (1) There is a vector $w^1\in\mathbb{R}^l$ such that the primal slackness vector \begin{equation} \label{prifea1} x=d+M^Tv+B^Tw^1\end{equation} is feasible for (\ref{primaljihe4}) if and only if $v\in \varTheta_P$. Moreover, when $d\in \operatorname{int}(K)$, $x$ is strictly feasible if and only if $v\in \operatorname{int}(\varTheta_P)$;

(2) There is a vector $w^2\in\mathbb{R}^m$ such that dual slackness vector \begin{equation} \label{prifea2} y=c+M^Tu+A^Tw^2\end{equation} is feasible for (\ref{primaljihe3}) if and only if $u\in \varTheta_D$. Moreover, when $c\in \operatorname{int}(K^*)$, $y$ is strictly feasible if and only if $u\in \operatorname{int}( \varTheta_D)$. \upshape

\prof The first part of the first claim follows from Assumption 1. It is obvious that $v\in \operatorname{int}( \varTheta_P)$ if $x$ is strictly feasible for (\ref{primaljihe4}). Conversely, if $d\in \operatorname{int}(K)$, then, for every $v\in \operatorname{int}( \varTheta_P)$, $x$ is strictly feasible. The proof is completed. $\square$

It should be noted that $d\in \operatorname{int}(K)$ implies that (\ref{primaljihe1}) is strictly feasible. Such a condition is necessary for the first argument of Corollary \ref{xinzd1}. We illustrate this with a semidefinite system, with $\mathbb{R}^{\frac{n(n+1)}{2}}\simeq\mathbb{S}^n$ as the set of order $n$ symmetric matrices and $K=K^*=\mathbb{S}^n_+$ as the set of order $n$ symmetric positive semidefinite matrices. The inner product of $\mathbb{S}^n$ is called $c \bullet d= \langle c,d\rangle=tr(cd)$. Note that we denote the elements of $\mathbb{S}^n$ in small letters. The row vectors of the matrix $A$ are denoted by $a^1,a^2,\cdots,a^m$. For any $x\in\mathbb{S}^n$, $x\succeq 0$ means $x\in\mathbb{S}^n_+$. To avoid complicated calculations, only Assumption 1 holds in all examples.

\exam \label{feasiblef1} Consider the following parametric SDP problem
\[ \begin{array}{cl} \min\limits_{x\in\mathbb{S}^3_+} & (c+m^1u)\bullet x \\ s.t. & a^i\bullet x = 0, \\ & a^2 \bullet x = 1, \end{array} \]
where $c=0\in\mathbb{R}^{3\times 3}$,
\[m^1=\left( \begin{array}{ccc} 1 &0&0 \\ 0&0&-0.5\\ 0&-0.5&0 \end{array} \right),\ a^1= \left( \begin{array}{ccc} 0 &0&0 \\ 0&1 &0\\ 0&0&0 \end{array} \right), \ a^2= \left( \begin{array}{ccc} 1 &0&0 \\ 0&0&1\\ 0&1&0 \end{array} \right). \] Then
\[ c+m^1u+a^1w_1+a^2w_2= \left( \begin{array}{ccc} u+w_2 &0&0 \\ 0&w_1&-0.5u+w_2\\ 0&-0.5u+w_2&0 \end{array} \right) \succeq 0\] for some $w_1,w_2\in\mathbb{R}$ if and only if
\[ u+w_2\geq 0, \quad w_1\geq 0,\quad -0.5u+w_2=0. \] From the first inequality and the last equality, we have $1.5u\geq0$, which means that $\varTheta_D =[0,+\infty)$. It has at least one interior point.

There is no feasible interior in the primal program, and for any $u\in\varTheta_D$, the nonstandard dual program
\[ \begin{array}{cl} \max\limits_{y\in\mathbb{S}^3_+} & d\bullet (c+m^1u-y) \\ s.t. & b^i\bullet y = 0, \quad i=1,2,3,\\ & m^1\bullet y =1.5 u \end{array} \] has no feasible interior, where
\[d = \left( \begin{array}{ccc} 1 &0&0 \\ 0&0&0\\ 0&0&0 \end{array} \right),b^1=\left( \begin{array}{ccc} 0 &1&0 \\ 1&0&0\\ 0&0&0 \end{array} \right), b^2=\left( \begin{array}{ccc} 0 &0&1 \\ 0&0&0\\ 1&0&0 \end{array} \right), b^3= \left( \begin{array}{ccc} 0 &0&0 \\ 0&0&0 \\ 0&0&1 \end{array} \right). \] Here $d$ is feasible, and can be replaced by other arbitrary primal feasible solutions. We ignore this description in the following statement.

Let $\mathscr{X}$ and $\mathscr{Y}$ denote the feasible sets in (\ref{primaljihe1}) and (\ref{primaljihe2}), respectively. These two sets do not depend on the given parametric vectors, $u$ and $v$. We assume that they are non-empty throughout this study. $\varTheta_P$ and $\varTheta_D$ are projections of $\mathscr{X}$ and $\mathscr{Y}$, respectively. In our terminology, $\varTheta_P$ and $\varTheta_D$ denote {\it primal and dual conic representable sets}, respectively.

Without confusion, the optimal solutions of (\ref{primaljihe1}) and (\ref{primaljihe2}) are denoted by $x^*(u)$ and $y^*(v)$, respectively, if they exist. That is,
\begin{eqnarray*} x^*(u) &\in& \mathscr{X}^*(u)= \arg\min \{\langle c+M^Tu,x\rangle | x\in \mathscr{X}\}, \\ y^*(v) & \in & \mathscr{Y}^*(v)= \arg\min \{\langle d+M^Tv,y\rangle | y\in \mathscr{Y}\}. \end{eqnarray*} Analogously, the optimal solutions of (\ref{primaljihe3}) and (\ref{primaljihe4}) are denoted as $\bar{y}^*(u)$ and $\bar{x}^*(v)$, respectively, if they exist. Or equivalently,
\begin{eqnarray*} \bar{y}^*(u) &\in& \bar{\mathscr{Y}}^*(u)= \arg\min \{\langle d,y\rangle | y \in \bar{\mathscr{Y}}(u)\}, \\ \bar{x}^*(v) & \in & \bar{\mathscr{X}}^*(v)= \arg\min \{\langle c,x\rangle | x \in \bar{\mathscr{X}}(v) \}, \end{eqnarray*}
where $\bar{\mathscr{X}}(v)$ and $\bar{\mathscr{Y}}(u)$ denote the feasible sets in (\ref{primaljihe4}) and (\ref{primaljihe3}), respectively.

\section{A modern duality} \label{t31}
This section presents a novel duality in CLO, namely the almost duality between (\ref{primaljihe1}) and (\ref{primaljihe2}) or (\ref{primaljihe4}) and (\ref{primaljihe3}).

\subsection{Set-valued mappings of the projections} \label{jzyss}
Our main aim in this subsection is to study how $\mathscr{X}$ and $\mathscr{Y}$ are projected onto the affine sets $\{x\in\mathbb{R}^q| Mx=Md+v\}$ and $\{y\in\mathbb{R}^q| My=Mc+u\}$, respectively. Let us begin with two rather technical lemmas that emerge from projection behaviors.

\lem \label{jshuxi1} (1) If $Mx = Md+v$ for $x\in \mathscr{X}$, then \begin{equation} \label{bdgh1} \langle M^Tu,x\rangle=\langle M^Tu,d+M^Tv\rangle.\end{equation}

(2) If $My = Mc+u$ for $y\in \mathscr{Y}$, then \begin{equation} \label{bdgh2} \langle M^Tv,y\rangle=\langle M^Tv,c+M^Tu\rangle.\end{equation} \upshape

\prof If $Mx = Md+v$ for some $x\in \mathscr{X}$, then
\begin{eqnarray*} \langle M^Tu,x\rangle = \langle u,Mx\rangle = \langle u,Md+v\rangle = \langle M^Tu,d+M^Tv\rangle, \end{eqnarray*} where $MM^T=I_r$ (Assumption 2) is used. The first claim holds. $\square$

Noted that Assumption 2 in Lemma \ref{jshuxi1} is unnecessary. We can replace $u$ (or $v$) with $(MM^T)^{-1}u$ (or $(MM^T)^{-1}v$), if this assumption does not hold.

\lem \label{jshuxi201} (1) For any $x\in \bar{\mathscr{X}}(v)$, the sum of the objective value of (\ref{primaljihe1}) and (\ref{primaljihe4}) is equal to $\langle c+M^Tu,d+M^Tv\rangle$, i.e.,
\begin{equation} \label{equalzeo1} \langle c+M^Tu, x \rangle + \langle c,d+M^Tv-x \rangle = \langle c+M^Tu,d+M^Tv\rangle. \end{equation}

(2) For any $y\in \bar{\mathscr{Y}}(u)$, the sum of the objective values of (\ref{primaljihe2}) and (\ref{primaljihe3}) is equal to $\langle c+M^Tu,d+M^Tv\rangle$, that is,
\begin{equation} \label{equalzeo2} \langle d+M^Tv, y \rangle + \langle d,c+M^Tu-y \rangle = \langle c+M^Tu,d+M^Tv\rangle. \end{equation} \upshape

\prof If $x\in \bar{\mathscr{X}}(v)$, then $x\in \mathscr{X}$ and $Mx=Md+v$. Furthermore, it follows from equality (\ref{bdgh1}) that equality (\ref{equalzeo1}) holds. The proof is completed. $\square$

The projection behaviors allow us to define two set-valued mappings as follows:
\begin{equation} \label{setvalue1} \Phi(u)= \left\{ M(x^*(u)-d)|x^*(u)\in \mathscr{X}^*(u) \right\}, \quad \forall u\in\varTheta_D \end{equation} and
\begin{equation} \label{setvalue2} \Psi(v)= \left \{ M(y^*(v)-c)|y^*(v)\in \mathscr{Y}^*(v) \right \}, \quad \forall v\in\varTheta_P. \end{equation}
 Here the value of every mapping can be a set if the optimal solution is not unique corresponding to the fixed parametric vector. We refer to \cite{AF09,Roc14,RW09} for a detailed introduction, and to \cite{Ber971,Bor03,Hog73} for the applications in optimization.

\thm \label{maithm1} (1) If $(d,c)\in K\times \operatorname{int}(K^*)$, then for every $u\in\operatorname{int}\left(\varTheta_D\right)$, $\Phi(u)$ is well-defined, i.e., $\Phi(u)\ne \emptyset$;

(2) If $(d,c)\in \operatorname{int}(K)\times K^*$, then for every $v\in \operatorname{int} \left(\varTheta_P\right)$, $\Psi(v)$ is well-defined, that is, $\Psi(v)\ne \emptyset$. \upshape

\prof If $d\in K$, then (\ref{primaljihe1}) has a feasible solution $d$. If $c \in \operatorname{int}( K^*)$, then $\varTheta_D$ has at least one interior $\bar{u}=0$, that is, $\operatorname{int}\left(\varTheta_D\right)\ne \emptyset$. From Corollary \ref{xinzd1}, for every $u\in\operatorname{int}\left(\varTheta_D\right)$, (\ref{primaljihe3}) satisfies Slater's condition. By Theorem \ref{dualityt1}, (\ref{primaljihe1}) is solvable, that is, $\mathscr{X}^*(u)\ne \emptyset$. Then $\Phi(u)$ is well defined. $\square$

In LP, both $\Phi(u)$ on $\varTheta_D$ and $\Psi(v)$ on $\varTheta_P$ are always well-defined. However, in SDP they could be undefined at the boundary of the conic representable sets. Here, is a counterexample.

\exam \label{patex10} Consider the following parametric SDP problem
\[ \begin{array}{ll} \min\limits_{x\in\mathbb{S}^2_+} & (c+m^1u) \bullet x \\ s.t. &
 a^1 \bullet x=2, \end{array} \] where \[c=\left(\begin{array} {cc} 1 & 0 \\ 0 & 0 \end{array} \right), m^1=\left(\begin{array} {cc} 0 & 0 \\ 0 & 1 \end{array} \right), a^1=\left(\begin{array} {cc} 0 & 1 \\ 1 & 0 \end{array} \right).\] Then for some $w\in\mathbb{R}$
 \[c+m^1u+a^1w=\left(\begin{array} {cc} 1 & w \\ w & u \end{array} \right)\succeq 0\]
means that $u\geq 0$, i.e., $\varTheta_D=[0,+\infty)$. The nonstandard dual program is
\[ \begin{array}{ll} \max\limits_{y\in\mathbb{S}^2_+} & d \bullet (c+m^1u-y) \\ s.t. & b^1\bullet y=1,\\ &m^1\bullet y=u, \end{array} \]where
\[d=\left(\begin{array} {cc} 2 & 1 \\ 1 & 2 \end{array} \right), \quad b^1=\left(\begin{array} {cc} 1 & 0 \\ 0 & 0 \end{array} \right). \]

If $u>0$, the primal and dual programs have optimal solutions \[ x^*(u)=\left( \begin{array}{cc}\sqrt{u} &1\\ 1&\frac{1}{\sqrt{u}} \end{array}\right),\qquad \bar{y}^*(u)=\left( \begin{array}{cc}1 &-\sqrt{u}\\ -\sqrt{u}& u \end{array}\right), \] respectively, and there is no duality gap. However, if $u=0$, the primal program is not solvable, although the dual program has 0 maximum at \[y^*(0)=\left( \begin{array}{cc}1 &0\\ 0&0\end{array}\right)\] and there is no duality gap. Therefore, for any $u\in (0,+\infty)$,
\[\Phi(u) = \left(\begin{array} {cc} 0 & 0 \\ 0 & 1 \end{array} \right)\bullet \left( \begin{array}{cc}\sqrt{u} &1\\ 1&\frac{1}{\sqrt{u}} \end{array}\right)-\left(\begin{array} {cc} 0 & 0 \\ 0 & 1 \end{array} \right)\bullet \left( \begin{array}{cc}2 &1\\ 1& 2\end{array}\right)
=\frac{1}{\sqrt{u}}-2 \] is well-defined and $\Phi(0)$ is undefined. In addition, it is easy to verify that the $\Psi(v)=\frac{1}{(v+2)^2}$ for any $v=(-2,+\infty)$ and $\Psi(-2)$ is undefined.

\cor \label{mainthe1} (1) Suppose that $\mathscr{X}^*(u)\ne \emptyset$ for some $u\in\varTheta_D$. Then $x^*(u)\in \bar{\mathscr{X}}^*(v)$ if and only if $v\in\Phi(u)$;

(2) Suppose $\mathscr{Y}^*(v)\ne \emptyset$ for $v\in\varTheta_P$. Then $ y^*(v)\in \bar{\mathscr{Y}}^*(u)$ if and only if $u\in\Psi(v)$. \upshape

\prof Suppose that $\mathscr{X}^*(u)\ne \emptyset$ for some $u\in\varTheta_D$. If $v\in \Phi(u)$, $x^*(u)\in\mathscr{X}^*(u)$ exists such that $v =M(x^*(u)-d)$. However, if $x\in \bar{\mathscr{X}}(v)$, $x\in \mathscr{X}$. From Lemma \ref{jshuxi201}, we have
\[\langle c+M^Tu, x \rangle + \langle c,d+M^Tv-x \rangle = \langle c+M^Tu, x^*(u) \rangle + \langle c,d+M^Tv-x^*(u) \rangle.\] Therefore,
\[\langle c,d+M^Tv-x^*(u) \rangle- \langle c,d+M^Tv-x \rangle = \langle c+M^Tu, x \rangle -\langle c+M^Tu, x^*(u) \rangle \geq 0,\] where the last inequality holds because of $x^*(u)\in \mathscr{X}^*(u)$. Therefore, the left-hand expression is greater than or equal to zero, that is, $x^*(u)\in \bar{\mathscr{X}}^*(v)$.

Conversely, if $x^*(u)\in \bar{\mathscr{X}}^*(v)$, $x^*(u)\in \bar{\mathscr{X}}^*(v)\cap\mathscr{X}^*(u)$. Then, by the second constraint of (\ref{primaljihe4}), we obtain $v\in\Phi(u)$. $\square$

Corollary \ref{mainthe1} gives an interesting geometric interpretation of the lift-and-project procedure. We illustrate it with a pair of lift-and-project problems (\ref{primaljihe1}) and (\ref{primaljihe4}). On the one hand, $\varTheta_D$ can denote the objective perturbation set of the lifted problem (\ref{primaljihe1}), in which the hyperplane $H_u=\{x\in \mathbb{R}^q| \langle c+M^Tu,x\rangle = \langle c+M^Tu,x^*(u)\rangle \}$ supports $\mathscr{X}$ at $x^*(u)$. On the other hand, $\varTheta_P$ denotes the right side hand perturbation set of (\ref{primaljihe4}), in which the affine hyperplane $C_v= \{x\in \mathbb{R}^q| Mx=Md+v\}$ cuts $\mathscr{X}$ at $x^*(u)$ if $v\in\Phi(u)$. Therefore, the set-valued mapping $\Phi(u)$ from $\varTheta_D$ to $\varTheta_P$ establishes the relationship between two different perturbations. For a given vector $u\in \varTheta_D$, if $H_u\cap \mathscr{X}$ is a set, then the cutting hyperplane $C_v$ is not unique, which implies that $\Phi(u)$ is a set. And if $x^*(u)$ is a vertex of $\mathscr{X}$, then the supporting hyperplane $H_u$ at $x^*(u)$ is not unique, which implies that there is a set $\mathcal{U}\subset \varTheta_D$ such that $\Phi(u)$ takes the same image for all $u\in \mathcal{U}$. When neither of the preceding cases occurs, there is only one the supporting hyperplane $H_u$ supporting the set $\mathscr{X}$, and there is only one hyperplane $C_v$ cutting the set $\mathscr{X}$. All of these cases will be examined in more detail later in the next section.

\subsection{Duality properties} \label{thro23}
In the classical Lagrangian dual approach, the complement slackness property is in demand in the primal and dual optimization pair such that the zero duality gap holds. However, for almost all primal-dual problems, their duality is slightly different from that of classical duality.

\cor \label{sumtwop21} {\bf (Weak duality and strong duality).} Suppose that $(d,c)\in K\times K^*$, $u\in \varTheta_D$ and $v\in \Phi(u)$ or $v\in \varTheta_P$ and $u\in \Psi(v)$.

(1) For any $(x,y)\in\bar{\mathscr{X}}(v)\times \bar{\mathscr{Y}}(u)$, one has
\begin{equation} \label{weadu1} \langle c, d+M^Tv-x\rangle + \langle d,c+M^Tu-y\rangle \leq \langle c+M^Tu, d+M^Tv\rangle.\end{equation} The equality holds if and only if $(x,y)\in\bar{\mathscr{X}}^*(v)\times \bar{\mathscr{Y}}^*(u)$.

(2) For any $(x,y)\in\mathscr{X}\times \mathscr{Y}$, one has
\begin{equation} \label{weadu2} \langle c+M^Tu,x\rangle + \langle d+M^Tv,y\rangle \geq \langle c+M^Tu, d+M^Tv\rangle. \end{equation} The equality holds if and only if $(x,y)\in\mathscr{X}^*(u)\times \mathscr{Y}^*(v)$. \upshape

\prof This result follows from Lemma \ref{jshuxi201} and Corollaries \ref{weakdual1} and \ref{mainthe1}. $\square$

The following theorem covers the complementary slackness properties (\ref{slack1}) and (\ref{slack2}).

\thm \label{interior2yy} {\bf (mpKKT property).} Let $u\in \varTheta_D$ and $v\in\varTheta_P$ be arbitrary. Then, there is a pair of vectors $(\bar{x},\bar{y})\in\mathbb{R}^q \times \mathbb{R}^q$ such that the mpKKT conditions hold
\begin{subequations}
\begin{align} \label{kktcond1} &A \bar{x} =b,\quad M \bar{x} =Md+v, \quad \bar{x} \in K, \\ \label{kktcond2} &
B\bar{y} = a,\quad M\bar{y} =Mc+u, \quad \bar{y} \in K^*,\\ \label{kktcond3} & \langle \bar{x},\bar{y} \rangle=0,
\end{align}
\end{subequations} if and only if both $v\in\Phi(u)$ and $u\in \Psi(v)$ hold. Furthermore, $\bar{x}\in \mathscr{X}^*(u)\cap \bar{\mathscr{X}}^*(v)$ and $\bar{y}\in \mathscr{Y}^*(v)\cap \bar{\mathscr{Y}}^*(u)$. \upshape

\prof We first show that if the mpKKT conditions (\ref{kktcond1})-(\ref{kktcond3}) hold, then both $v\in\Phi(u)$ and $u\in \Psi(v)$ hold. Clearly, the above mpKKT conditions imply that \begin{equation} \label{addtj1} M\bar{y} =Mc+u\end{equation} and \begin{eqnarray*} &A \bar{x}=b,\quad M \bar{x} =Md+v, \quad \bar{x} \in K, \\ &
B\bar{y} = a, \quad \bar{y} \in K^*,\\ & \langle \bar{x},\bar{y} \rangle=0. \end{eqnarray*} From the KKT conditions (\ref{jbkkt1})-(\ref{jbkkt3}), the latter implies that $(\bar{x},\bar{y})$ is a pair of optimal solutions of (\ref{primaljihe4}) and (\ref{primaljihe2}). Therefore, $\bar{y}\in \mathscr{Y}^*(v)$, i.e., there is $y^*(v)\in \mathscr{Y}^*(v)$ such that $y^*(v)=\bar{y}$. Finally, equality (\ref{addtj1}) implies that $u=M(y^*(v)-c)\in\Psi(v)$. Analogously, one has $v\in\Phi(u)$.

The converse is easy to prove. If both $v\in\Phi(u)$ and $u\in \Psi(v)$ hold, then $\Phi(u)$ and $\Psi(v)$ are well defined, which means that almost primal-dual programs (\ref{primaljihe1}) and (\ref{primaljihe2}) have a pair of optimal solutions $(x^*(u),y^*(v))$. Furthermore, by Corollary \ref{mainthe1}, $(x^*(u),y^*(v))$ is also a pair of optimal solutions for (\ref{primaljihe4}) and (\ref{primaljihe3}). Thus, from the KKT conditions (\ref{jbkkt1})-(\ref{jbkkt3}), $(\bar{x},\bar{y})=(x^*(u),y^*(v))$ satisfies mpKKT conditions (\ref{kktcond1})-(\ref{kktcond3}).

The final claim follows from KKT conditions (\ref{jbkkt1})-(\ref{jbkkt3}). $\square$

Our duality inherits many useful properties of classical duality, although some are different from classical forms. Because it depends on the perturbation parameters, duality presents more dynamic features of a CLO problem. If two pairs of almost primal-dual programs are put together, then these features, which are well connected by set-valued mappings, have a very intuitive geometric explanation. For example, refer to the remarks at the end of the previous subsection. In short, duality embodies a dynamic process, whereas classical duality is a static result.

\subsection{Other related results} \label{otjg1}
The following result shows that the set-valued mapping $\Phi(u)$ iterates over every value in the set $\operatorname{int}(\varTheta_P)$ if $u$ goes through every value in the set $\varTheta_D$, although $\Phi(u)$ could be undefined for some boundary point of $\varTheta_D$.

\thm \label{interior1} Suppose that $(d,c)\in \operatorname{int}(K\times K^*)$, then
\begin{subequations}
\begin{align}\label{baohan1} \operatorname{int}(\varTheta_P) \subset \bigcup\limits_{u\in\varTheta_D} \Phi(u), \\ \label{baohan2} \operatorname{int}(\varTheta_D) \subset \bigcup\limits_{v\in\varTheta_P} \Psi(v).
\end{align}
\end{subequations} \upshape

\prof If $v\in \operatorname{int}(\varTheta_P)$, then (\ref{primaljihe4}) and (\ref{primaljihe2}) satisfy the Slater's conditions. From Theorem \ref{dualityt1}, there are $ \bar{x}^*(v) \in \bar{\mathscr{X}}^*(v) $ and $y^*(v)\in \mathscr{Y}^*(v)$ such that
\[\langle c,d+M^Tv-\bar{x}^*(v)\rangle=\bar{p}^*(v)=d^*(v)=\langle d+M^Tv,y^*(v)\rangle. \] From Corollary \ref{mainthe1}, $u\in \varTheta_D$ exists such that $u\in \Psi(v)$. It follows from (\ref{equalzeo1}), (\ref{weadu1}) and (\ref{weadu2}) that
 \[p^*(u)= \langle c+M^Tu,d+M^Tv\rangle-d^*(v) = \langle c+M^Tu,d+M^Tv\rangle-\bar{p}^*(v)= \langle c+M^Tu,\bar{x}^*(v)\rangle. \] Because it is feasible for (\ref{primaljihe1}), $\bar{x}^*(v)$ is optimal for (\ref{primaljihe1}), that is, $\bar{x}^*(v)\in \mathscr{X}^*(u)$. Or equivalently, $v\in\Phi(u)$. The proof is completed. $\square$

Sometimes, we require a version of Theorem \ref{interior1}. This version shows that set-valued mapping $\Psi(v)$ is almost the inverse of set-valued mapping $\Phi(u)$.

\cor \label{interior2} Suppose that $(d,c)\in \operatorname{int}(K\times K^*)$.

(1) If $u\in \Psi(v)$ for some $v\in \operatorname{int}(\varTheta_P)$, then $v\in \Phi(u)$.

(2) If there is $v\in \Phi(u)$ for some $u\in \operatorname{int}(\varTheta_D)$, then $u\in \Psi(v)$. \upshape

 \cor \label{xinzd3} Suppose that $(d,c)\in \operatorname{int}(K\times K^*)$.

(1) For every $u\in \operatorname{int} \left(\varTheta_D\right)$, $\Phi(u)$ is a closed convex set. This can be described as \begin{equation} \label{biaoshi1} \Phi(u)=\{v| \operatorname{\exists \ (\bar{x},\bar{y})\ s.t. \ the\ mpKKT\ conditions} \ (\ref{kktcond1})-(\ref{kktcond3})\ \operatorname{holds}\}; \end{equation}

(2) For every $v\in \operatorname{int} \left(\varTheta_P\right)$, $\Psi(v)$ is a closed convex set. This can be identified by \begin{equation} \label{biaoshi2} \Psi(v) = \{u | \operatorname{\exists \ (\bar{x},\bar{y})\ s.t. \ the\ mpKKT\ conditions} \ (\ref{kktcond1})-(\ref{kktcond3})\ \operatorname{holds}\} \end{equation} \upshape
\indent \prof Firstly, the equalities (\ref{biaoshi1}) and (\ref{biaoshi2}) follow from Theorem \ref{interior2yy} and Corollary \ref{interior2}.

Now we prove that $\Phi(u)$ for every $u\in \operatorname{int} \left(\varTheta_D\right)$ is a closed convex set. If $v^1,v^2\in \Phi(u)$, then from (\ref{biaoshi1}) and (\ref{jbkkt1})-(\ref{jbkkt3}), there are $\bar{x}^1, \bar{x}^2$ and $\bar{y}^*$ such that
\[ \begin{array}{llll} A\bar{x}^1=b, & M \bar{x}^1=Mc+v^1, & \bar{x}^1\in K, \\ B\bar{y}=a, & M\bar{y} = Mc + u, & \bar{y}\in K^*, \\ \langle \bar{x}^1,\bar{y}\rangle=0 \end{array}\] and
\[ \begin{array}{llll} A\bar{x}^2=b, & M \bar{x}^2=Mc+v^2, & \bar{x}^2\in K, \\ B\bar{y}=a, & M\bar{y} = Mc + u, & \bar{y}\in K^*, \\ \langle \bar{x}^2,\bar{y}\rangle=0. \end{array}\]
Therefore, for any $\alpha\in [0,1]$, $\bar{x}_{\alpha}=\alpha \bar{x}^1+(1-\alpha)\bar{x}^2$ satisfies
\[ \begin{array}{llll} A\bar{x}_{\alpha}=b, & M \bar{x}_{\alpha}=Mc+\alpha v^1+(1-\alpha)v^2, & \bar{x}_{\alpha}\in K, \\ B\bar{y}=a, & M\bar{y} = Mc + u, & \bar{y}\in K^*, \\ \langle \bar{x}_{\alpha},\bar{y}\rangle=0. \end{array}\] Applying Theorem \ref{interior2yy} again, we obtain $v=\alpha v^1+(1-\alpha)v^2\in \Phi(u)$. That is, $\Phi(u)$ is convex.

Finally, the intersection of the supporting hyperplane is obtained
\[ H =\{x\in\mathbb{R}^q| \langle c+M^Tu,x\rangle = \langle c+M^Tu,x^*(u)\rangle \}\] and $ X=\{x\in K|Ax=b\}$ is closed, where $x^*(u)$ equals $\bar{x}^1$ or $\bar{x}^2$. This means that $\Phi(u)$ is closed. The proof is finished. $\square$

\cor \label{induct1} Let $u\in \varTheta_D$ and $v\in \varTheta_P$ be both arbitrary.

(1) Suppose the $(d,c)\in K\times \operatorname{int}(K^*)$. (\ref{primaljihe4}) is solvable if and only if it is $v\in\Phi(u)$;

(2) Suppose the $(d,c)\in \operatorname{int}(K)\times K^*$. (\ref{primaljihe3}) is solvable if and only if it is $u\in\Phi(v)$. \upshape

\prof Let us prove the first claim. By Theorem \ref{dualityt1}, (\ref{primaljihe1}) and (\ref{primaljihe4}) are solvable, in which the cutting hyperplane $C_v=\{x\in\mathbb{R}^q|Mx=Md+v\}$ must intersect the feasible set $ X=\{x\in K|Ax=b\}$. Then, from Corollary \ref{mainthe1}, (\ref{primaljihe4}) is solvable if and only if $v\in\Phi(u)$. $\square$

To conclude this subsection, we offer a new proof of Theorem \ref{dualityt1} using Corollary \ref{induct1}.

Proof of Theorem \ref{dualityt1}. Let us assume that (\ref{primaljihe30}) is bounded and strictly feasible, that is, $(d,c)=K\times \operatorname{int}(K^*)$. We prove the first claim through mathematical induction for $k=q-m$.

Initial step: $k=0$. Here, the feasible set in (\ref{primaljihe10}) is a singleton set. It has a minimum and there is no duality gap. The first claim holds.

Inductive step: Assume that (\ref{primaljihe4}) is solvable for some $v\in\varTheta_P$. From Corollaries \ref{mainthe1} and \ref{induct1},
(\ref{primaljihe1}) for all $u\in \operatorname{int}(\varTheta_D)$ is solvable, and there is no duality gap. Because $0\in \operatorname{int}(\varTheta_D)$, (\ref{primaljihe10}) is solvable, there is no duality gap for $k:=k-l$. The first claim then holds by induction. $\square$

If $L=N(A)$, or equivalently $L=R(M^T)\oplus R(B^T)$ (Ref \cite{HJ89}), then we can rewrite (\ref{primaljihe10}) in the geometric form
\begin{equation} \label{geome1} \min \langle c,x\rangle,\quad s.t. \ x \in d+L,\quad x\in K. \end{equation} Correspondingly, its dual is in the geometric form
\begin{equation} \label{geome2} \min \langle d,y\rangle,\quad s.t. \ y \in c+L^{\bot},\quad y\in K^*,\end{equation} where $L^{\bot}$ denotes the orthogonal complement of $L$ in $\mathbb{R}^q$, e.g., see Nesterov and Nemirovski \cite{NN94}, and Todd \cite{Tod09}. Under Assumption 1, $L^{\bot}=R(A^T)=N(B)\cap N(M)$ implies that the aforementioned dual program (\ref{geome2}) is the same as (\ref{primaljihe30}). Of course, the choice of matrices $B$ and $M$ for (\ref{primaljihe30}) is not unique when the matrix $A$ is fixed.

In the proof of Theorem \ref{dualityt1}, the row vectors of $M,A,B$ denote the row vectors of the changed, old primal and new dual constraint matrices, respectively. In the lift-and-project procedure, $L^\bot$ is lifted from $N(B)\cap N(M)$ to $N(B)$, and $L$ projects from $N(A)$ to $N(A)\cap N(M)$. In other words, $R(M^T)$ is not only orthogonal to the old primal space $R(A^T)$, but also to the new dual space $R(B^T)$. This motivated us to propose Assumption 1.

\section{Invariancy sets and and illustrative examples}
The concept of optimal partition was introduced for parametric LP in \cite{AD92,JR03}, in which the given optimal basic partition is invariant. Later, this concept was discussed in \cite{GY93,Wri08} for quadratic programming problems, \cite{Ber97} for linear complementarity problems and \cite{GS99,MT20} for SDPs. In this section we use the set-valued mappings described in the previous section to define the optimal partition of conic representable sets.

\defn \label{djkyi1} Let $\mathcal{V}$ be a simply connected subset of $\varTheta_P$. Then

(1) $\mathcal{V}$ is called a {\it linearity set} if $\Psi(v^1)=\Psi(v^2)$ for all $v^1,v^2\in\mathcal{V}$, in which $\mathcal{V}$ and $\Psi(\mathcal{V})$ are at least not a singleton set.

(2) $\mathcal{V}$ is called a {\it nonlinearity set} if $\mathcal{V}$ is not a singleton set, and for any $v^2\ne v^1\in\mathcal{V}$, $\Psi(v^2)\ne \Psi(v^1)$, and $\Psi(v^1)$ is a singleton set.

\noindent Both a linearity set and a nonlinearity set are called an {\it invariancy set}. For the dual conic representable set $\varTheta_D$, the definitions of the invariant set and the linear/nonlinearity set are similar.

From Definition \ref{djkyi1}, a conic representable set mainly includes two types of invariancy sets: linearity and nonlinearity sets. For every linearity set, the primal (or dual) optimal solution remains unchanged if the $\Psi(\mathcal{V})$ (or $\Phi(\mathcal{U})$) is a singleton set. In LP, there are only linearity sets; however, in SDP, there could be a nonlinearity set that was carefully studied by Mohammad-Nezhad and Terlaky \cite{MT20}. As indicated in \cite{GG08,MT20} and also demonstrated by Example \ref{patex10}, the optimal partition for $\varTheta_P$ may vary with the parameter $v$ on a subinterval of $(-2,+\infty)$; and the optimal partition for $\varTheta_D$ may vary with the parameter $u$ on a subinterval of $[0,+\infty)$. It is easy to see that $\Phi(u)$ on $(0,+\infty)$ and $\Psi(v)$ on $(-2,+\infty)$ are a pair of inverse functions; and by this, $(-2,+\infty)$ and $(0,+\infty)$ the nonlinearity intervals of $\varTheta_P$ and $\varTheta_D$, respectively.

\defn Let $\mathcal{V}$ be an invariancy set of $\varTheta_P$. $\mathcal{V}$ is called a transition face if $\operatorname{dim}(\mathcal{V})<r$. In particular, if $\operatorname{dim}(\mathcal{V})=0$, then $\mathcal{V}$ is called a transition point; and if $\operatorname{dim}(\mathcal{V})=1$, $\mathcal{V}$ is called a transition line, etc. \upshape

We refer to an invariancy set as a nontrivial invariancy set if it is not a transition face. By contrast, a transition face is called a trivial invariancy set. Following we state the main result of this section.

\thm \label{zyjg1} Any two different invariancy regions of a conic representable set do not intersect. \upshape

\prof This result follows immediately from Definition \ref{djkyi1}. $\square$

From Corollary \ref{mainthe1} and Theorem \ref{zyjg1}, we can obtain the invariant set decomposition of a conic representable set. To build intuition, we present the following two examples: one of them is an LP problem and an SDP problem.

\exam\label{xxgh1} Taking
\[\begin{array}{ccl} c&=& (-1,-1,0,0,0)^T,\\
A&=& \left(\begin{array}{ccccc} 1&1&1&0&0\\ 0&1&0&1&0 \\ 1&0&0&0&1 \end{array} \right), \\
M&=& (-0.5,0.5,0,-0.5,0.5),\\
b&=&(3,2,2.5)^T, \end{array}\]
the following parametric LP program
\[\begin{array}{cl} \min\limits_{x\in\mathbb{R}^5} & (c+M^Tu)^Tx\\ s.t. &Ax=b,\quad x\geq 0 \end{array} \]
has a feasible solution $d=(1,1,1,1,1.5)^T$. The non-standard dual program can be expressed as follows:
\[\begin{array}{cl} \max\limits_{y\in\mathbb{R}^5} & d^T(c+M^Tu-y) \\ s.t. &By=-2, My=u, y\geq 0, \end{array} \] where $B=(1,1,-2,-1,-1)$. If the slackness variables $x_3,x_4,x_5$ are omitted, the primal feasible set reduces the convex pentagon
\[\{x\in\mathbb{R}^2| x_1+x_2\leq 3, 0\leq x_1\leq 2.5, 0\leq x_2\leq 2\}\] in a two-dimensional plane. This convex pentagon has five vertices
\[(0,0)^T, (0,2)^T, (2.5,0)^T, (2.5,0.5)^T, (1,2)^T. \]
For $ 0\leq u \leq 1$, the optimal pair $(x^*(u),\bar{y}^*(u))$ is as follows
\[x^*(u)=(2.5,0.5,0,1.5,0)^T,\qquad \bar{y}^*(u)=(0,0,1-u,0,2u)^T. \] Geometrically, the trajectory of $\bar{y}^*(u)$ in the interval $[0,1]$ is an edge of the polyhedral
\[\{y\in\mathbb{R}^5| y_1+y_2-2y_3 -y_4-y_5 =-2, y_1,y_2,y_3,y_4,y_5\geq 0\}, \] where the edge connects two vertices: $y^1=(0,0,1,0,0)^T$ and $y^2=(0,0,0,0,2)^T$.
Then for every $u\in (0,1)$, $\Phi(u)=-2$. If $v=Mx^*(u)-Md=-2$, $\Psi(v)=[0,1]$ is an interval. When $u$ is equal to either $0$ or $1$, $\Phi(u)$ is also an interval.
The following parallel table lists the values of the set-valued mappings $\Phi(u)$ and $\Psi(v)$.
 \[\begin{array}{c} Table \ 1 \quad The \ set-valued\ mappings \ in \ Example\ \ref{xxgh1} \\ \begin{array}{|c|c|c|c|} \hline \bar{x}^*(v) & v & u & \bar{y}^*(u) \\ \hline (0,2,1,0,2.5)^T &2& (-\infty,-1) & (-1-u,0,0,1-u,0)^T \\ \hline (2-v,2,v-1,0,0.5+v)^T& (1,2)
 & -1 & (0,0,0,2,0)^T
 \\ \hline (1,2,0,0,1.5) & 1 &(-1,0)& (0,0,1+u,-2u,0)^T \\ \hline 0.5(3-v,3+v,0,1-v,2+v) & (-2,1) & 0 & (0,0,1,0,0)^T
 \\ \hline (2.5,0.5,0,1.5,0)^T & -2 & (0,1) & (0,0,1-u,0,2u)^T \\ \hline (2.5,v+2.5,-2-v,-0.5-v,0)^T & (-2.5,-2) & 1 & (0,0,0,0,2)^T \\ \hline
 (2.5,0,0.5,2,0)^T & -2.5 & (1,+\infty) &
 (0,0,u-1,0,u+1)\\ \hline \bar{x}^*(v) & \Phi(u) & \Psi(v) & \bar{y}^*(u) \\ \hline\end{array} \end{array}\]

The following observations can be made from Table 1.

1. The primal linear representable set $\varTheta_P$ is equivalent to $[-2.5,2]$. It contains three open invariancy intervals $(-2.5,-2),(-2,-1),(1,2)$ and four transition points $-2.5,-2,1,2$. The dual linear representable set, $\varTheta_D$, is equal to $(-\infty,+\infty)$. It contains four open invariancy intervals $(-\infty,-1),(-1,0),(0,1),(1,+\infty)$ and three transition points $-1,0,1$. For each invariancy interval, the trajectory of the optimal solution is an edge of the polyhedron connecting two adjacent vertices. For each transition point, the corresponding optimal solution is the vertex of the polyhedron.

2. All invariancy sets are linearity. These are either open invariancy intervals or endpoints of invariancy intervals, in which the endpoints are transition points. The image of the set-valued mapping is a closed interval if and only if the primage is a transition point, and it is a transition point if and only if the primage is an open interval.

3. The intersection of the images of set-valued mapping at any two transition points is either empty or a transition point associated with its dual representable set.

4. Either $\Phi((-\infty,-1))=2$ or $\Psi(2)=(-\infty,-1]$ implies that the primal objective function with perturbations takes the maxima and minima in the primal feasible set. In addition, both $\Phi((1,+\infty))=-2.5$ and $\Psi(-2.5)=[1,+\infty)$ imply that the perturbed objective function takes the maxima and minima in the dual feasible set.

In this example, the given optimal basic partition is invariant for single-parametric LPs, which is similar to the previous definition (e.g., see \cite{AD92}).

\begin{figure}[htbp]
\centering\includegraphics[width=3.5in]{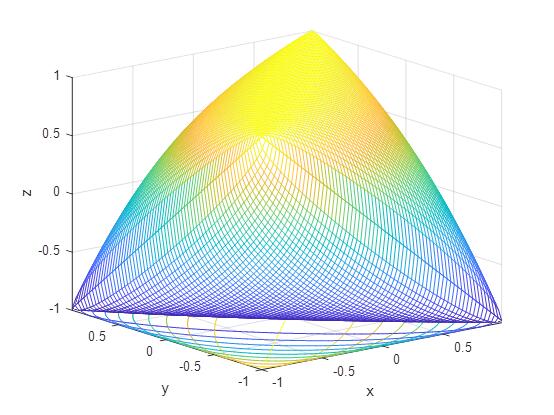}
\caption{The primal feasible set is a 3-elliptope}\label{fig:2}
\end{figure}

\exam \label{patex2} Take
\begin{eqnarray*} c=\left(\begin{array} {ccc} 0& 0 & 0 \\ 0 & 0 &-1 \\ 0 &-1 & 0 \end{array} \right), & m^1=\left(\begin{array} {ccc} 0& 1 & 0\\ 1 & 0 &0 \\ 0 &0 & 0 \end{array} \right), & m^2=\left(\begin{array} {ccc} 0& 0 & -1 \\ 0 & 0 &0 \\ -1 & 0& 0 \end{array} \right), \\
a^1 =\left(\begin{array} {ccc} 1 & 0 & 0 \\ 0& 0 & 0 \\ 0&0&0 \end{array} \right), & a^2=\left(\begin{array} {ccc} 0 & 0 & 0 \\ 0& 1 & 0 \\ 0&0&0 \end{array} \right), & a^3=\left(\begin{array} {ccc} 0 & 0 & 0 \\ 0& 0 & 0 \\ 0&0&1 \end{array} \right) \end{eqnarray*} and
\[ d=\left(\begin{array} {ccc} 1& 0 & 0 \\ 0 & 1 &0 \\ 0 &0 & 1 \end{array} \right), b^1=\left(\begin{array} {ccc} 0& 0 & 0 \\ 0 & 0 &1 \\ 0 &1 & 0 \end{array} \right). \] Consider the following mpSDP pair
 \begin{equation} \label{mostl1} \begin{array}{ll} \min\limits_{x\in\mathbb{S}^3_+} & (c+m^1u_1+m^2u_2)\bullet x\\ s.t. & a^i\bullet x=1,\quad i=1,2,3 \end{array} \end{equation} and \begin{equation} \label{mostl2} \begin{array}{ll} \min\limits_{y\in\mathbb{S}^3_+} & (d+m^1v_1+m^2v_2)\bullet y\\ s.t. & b^1\bullet y=-2. \end{array} \end{equation}
The primal feasible set is a 3-ellitope whose image is shown in Figure \ref{fig:2}.

The nonstandard dual problem of (\ref{mostl1}) is
\[\begin{array}{ll} \max\limits_{y\in\mathbb{S}^3_+} & d\bullet(c+m^1u_1+m^2u_2-y)\\ s.t. & b^1\bullet y=-2, \\ &m^1\bullet y =2u_1, \\ & m^2\bullet y =2u_2. \end{array}\]
The optimal solution pair $(x^*(u),\bar{y}^*(u))$ has six indeterminate entries:
\[ x^*(u) = \left(\begin{array} {ccc} 1 & x_{12} & x_{13} \\ x_{12} & 1 & x_{23} \\ x_{13}& x_{23} & 1\end{array} \right),\qquad
\bar{y}^*(u)=\left(\begin{array} {ccc} \bar{y}_{11} & u_1 & -u_2 \\ u_1 & \bar{y}_{22} & -1 \\ -u_2 &-1 & \bar{y}_{33} \end{array} \right), \]
in which we are interested in the six indeterminate entries, $x_{11},x_{12},x_{23}$ and $\bar{y}_{11},\bar{y}_{22},\bar{y}_{33}$, as a function of $(u_1,u_2)\in\mathbb{R}^2$.
Complex calculations are presented in Appendix A.

In this example, the dual conic representable set $\varTheta_D$ is the entire two-dimensional space and the primal
conic representable set $\varTheta_P$ is rectangular, given by \[ \varTheta_P=\{(v_1,v_2)^T| -1\leq v_1\leq 1,\ -1\leq v_2\leq 1\}.\]
For the set $\varTheta_D$, two open sets contained in $\varTheta_D^0$ are nonlinearity sets, and four angular domains $\varTheta_D^1,\varTheta_D^2,\varTheta_D^3-\{(0,0)^T\},\varTheta_D^4-\{(0,0)^T\}$ are nontrivial invariancy regions, the origin of which is only one transition point that is not associated with vertices. For the set $\varTheta_P$, the two open triangles in the two-dimensional plane are given by:
\[ \varTheta_P^0=\{(v_1,-v_2)^T| -1< v_1\ne v_2 < 1\}\] are nonlinearity sets and four angular points \[(1,1)^T, (-1,-1)^T, (1,-1)^T,(-1,1)^T\] are the transition points associated with the vertices. Diagonal rectangular region $\varTheta_P$ \[\varTheta_P^l=\{(v_1,-v_1)^T| -1< v_1< 1\}\] is a transition line of $\varTheta_P$, and its closure is a line segment that joins the two diagonal angular points $(1,-1)^T$ and $(-1,1)^T$. However, for any
\[v\in \partial (\varTheta_P)-\{(1,1)^T, (-1,-1)^T, (1,-1)^T,(-1,1)^T\}, \] one has $\Psi(v)=\emptyset$. A comparison of the set-valued mappings is shown in Table 2,
 \[\begin{array}{c} Table \ 2 \quad The \ set-valued\ mappings \ in \ Example \ \ref{patex2} \\ \begin{array}{||c||c|c|c|c|c||c||c||} \hline \hline v & (1,1)^T & (-1,-1)^T & (1,-1)^T & \varTheta_P^l &(-1,1)^T &\varTheta_P^0 & \Phi(u) \\ \hline u & \varTheta_D^1 & \varTheta_D^2 & \varTheta_D^3 & (0,0)^T = \varTheta_D^3\cap \varTheta_D^4 & \varTheta_D^4 & \varTheta_D^0 & \Psi(v) \\ \hline \hline\end{array} \end{array}\] where either $\varTheta_P^0$ or $ \varTheta_D^0$ contains two nontrivial invariancy sets; $\Phi$ and $\Psi$ over nonlinearity sets are given by formulas (\ref{setvalue11}) and (\ref{setvalue12}), respectively.

\section{Multiparametric analysis}
In this section, we present some applications of our duality properties to the multiparametric analysis of CLOs.

\subsection{Recession directions} \label{recess1}
The following corollary provides a theoretical foundation for the shadow vertex algorithm, see, e.g., \cite{Bor87,GS55}.

\cor \label{sensiti2} Suppose that $(d,c)\in \operatorname{int}(K\times K^*)$.

(1) If $h\in 0^+(\varTheta_D)$ is in a linearity set, then there is a vector $v^0\in \operatorname{\partial}(\varTheta_P)$ such that $h\in 0^+(\Psi(v^0))$;

(2) If $ h\in 0^+(\varTheta_P)$ is in a linearity set, then there is a vector $u^0\in \operatorname{\partial}(\varTheta_D)$ such that $h\in 0^+(\Phi(u^0))$. \upshape

\prof If $h\in 0^+(\varTheta_D)$, then (\ref{primaljihe3}) with $u=u^0+\lambda h$ is strictly feasible for any $u^0\in \varTheta_D$ and $\lambda>0$. Applying Theorem \ref{dualityt1}, (\ref{primaljihe1}) is solvable for $u=u^0+\lambda h$. If $h$ is in a linearity set, then there is a sufficiently large positive number $\Lambda$ such that for any $\lambda>\Lambda$, (\ref{primaljihe1}) with $u=u^0+\lambda h$ reaches a minimum at the same point $x^*(h)$ for some $u^0$. Then, (\ref{primaljihe4}) with $v^0= Mx^*(h)-Md\in \Phi(u^0+\lambda h)$ attain a minimum at $\bar{x}^*(v^0)=x^*(h)$. From Corollary \ref{interior2}, for any $\lambda>\Lambda$, $u=u^0+\lambda h\in \Psi(v^0)$. From Corollary \ref{xinzd3}, we obtain $h\in 0^+(\Psi(v^0))$. Finally, from the optimality of $\bar{x}^*(v^0)$, we obtain $v^0\in \operatorname{\partial}(\varTheta_P)$. $\square$

In Example \ref{xxgh1}, the two invariancy intervals, $(-\infty,-1]$ and $[1,+\infty)$, of the dual conic representable set contain the recession directions. In Example \ref{patex2}, the four invariancy regions, $\varTheta_D^1, \varTheta_D^2, \varTheta_D^3,\varTheta_D^4$, of the dual conic representable set contain the recession directions. There is no recession direction in the respective primal conic representable sets, because they are bounded. Finally, in Corollary \ref{sensiti2}, it is necessary to assume that recession directions belong to linearity sets. For an example, see Example \ref{patex10}.

\subsection{Identification of the optimal partitions}
When the parametric vectors $u$ and $v$ are scalar, the nontrivial invariancy set reduces the invariancy interval, and the transition face reduces the transition point. At the endpoints of the invariancy interval, the optimal partitions change upon transitioning to the adjacent invariancy interval. The transition point for parametric LPs, as a separation point of invariancy intervals, was mentioned many times (for example, see \cite{GG08,GM08,RT05}). The concepts of the transition point and nonlinearity invariant set for parametric SDPs were formally defined by Mohammad-Nezhad and Terlaky \cite{MT20}.

\thm \label{sensiti3} (1) If $(d,c)\in K \times \operatorname{int}(K^*)$, then every nonlinearity region $\mathcal{U}$ of $\varTheta_D$ is open and $\Phi(u)$ is continuous in the region $\mathcal{U}$;

(2) If $(d,c)\in \operatorname{int}(K)\times K^*$, then every nonlinearity region $\mathcal{V}$ of $\varTheta_P$ is open and $\Phi(v)$ is continuous in the region $\mathcal{V}$. \upshape

\prof We first show that the nonlinearity region is open. We assume that for every $\bar{u}\in \mathcal{U}$, $\Phi(\bar{u})$ is a singleton set. Geometrically, supporting hyperplane
\[ H_{\bar{u}} = \{x\in \mathbb{R}^q| \langle c+M^T\bar{u},x\rangle= \langle c+M^T\bar{u},x^*(\bar{u})\rangle\}\]
is tangent to the primal feasible set $ X=\{x\in K| Ax=b\}$ at the unique point $x^*(\bar{u})$, which implies that there is a neighbourhood $U(\bar{u})$ such that for any $u\in U(\bar{u})$, $x^*(u)$ lies in the local smooth surface of $ X$. Subsequently, for any $u\in U(\bar{u})$, the supporting hyperplane \[ H_u= \{x\in \mathbb{R}^q| \langle c+M^Tu,x\rangle= \langle c+M^Tu,x^*(u)\rangle\}\] is tangent to the primal feasible set $ X$. Furthermore, for any $u\in U(\bar{u})$, $\Phi(u)$ is a singleton set, meanings that $\mathcal{U}$ is open.

We now show that $\Phi(u)$ is continuous in the region $\mathcal{U}$. If for every $u\in \mathcal{U}$, $\Phi(u)$ is a singleton set, then the set-valued map $\Phi(u)$ on $\mathcal{U}$ degrades into a single-valued map. Then for two different vectors $u^1\in\mathcal{U}$ and $ u^2\in\mathcal{U}$, $\Phi(u^1)\ne \Phi(u^2)$. By the connectivity of $\mathcal{U}$, if $\Gamma$ denotes a continuous curve connecting two different points $u^1\in \mathcal{U}$ and $u^2\in \mathcal{U}$, then the trajectory of the optimal solution $x^*(u)$ ($u\in\Gamma$) is a continuous curve along the boundary of the primal feasible set. Therefore, $\Phi(\Gamma)$ is a continuous curve connecting two points $v^1=\Phi(u^1)$ and $v^2=\Phi(u^2)$, which implies the continuity of $\Phi(u)$ over the set $\mathcal{U}$. The proof is finished. $\square$

Every set-valued mapping in a nonlinearity set reduces to an ordinary single-valued mapping. By Corollaries \ref{interior2} and \ref{sensiti3}, $\Psi(v)$ and $\Phi(u)$ together form a pair of reversible single-valued mappings (for example, see Example \ref{patex2}).

It should be noted that our definition of the nonlinearity set differs from that given by Mohammad-Nezhad and Terlaky \cite{MT20}. In our definition, a nonlinearity set can not contain a transition point; however, in \cite{MT20}, a nonlinearity set could contain a transition point. For example, if $u_1=u_2$ in Example \ref{patex2}, then the origin belongs to the nonlinearity interval \cite[Example 3.1]{MT20}. Recently, Hauenstein et al. \cite{Hau19} analyzed the continuity of optimal set mapping and showed that continuity may fail on a nonlinearity interval under the old definition. This result does not contradict Theorem \ref{sensiti3}.

A direct consequence of Theorem \ref{sensiti3} is as follows.

\cor \label{tplog2} Let $\mathcal{U}$ and $\mathcal{V}$ be invariancy sets of $\varTheta_D$ and $\varTheta_P$, respectively.

(1) Suppose the $(d,c)\in K \times \operatorname{int}(K^*)$. Then, $\mathcal{U}$ is nonlinearity if and only if $\Phi(\mathcal{U})$ is nonlinearity;

(2) Suppose the $(d,c)\in \operatorname{int}(K)\times K^*$. Then, $\mathcal{V}$ is nonlinearity if and only if $\Psi(\mathcal{V})$ is nonlinearity. \upshape

In Example \ref{patex2}, for the primal conic representable set $\varTheta_P$, four vertices are transition points and the two open sets consisting of interiors of $\varTheta_P$ separated by the diagonal $\{(v_1,-v_1)^T| v_1\in (-2,2)\}$ are nonlinearity sets, in which the diagonal is a transition line of $\varTheta_P$. For the dual conic representable set $\varTheta_D$, the origin is only one transition point and the regions contained in the two curved triangles in the first and third quadrants are two different nonlinearity sets.

\thm \label{tplog1} Suppose that $(d,c)\in \operatorname{int}(K \times K^*)$.

(1) If $\mathcal{U}$ is a nontrivial linearity set of $\varTheta_D$, then $\mathcal{U}$ is convex, and for all $u\in \operatorname{cl}(\mathcal{U}) $ and $v\in \Phi(\mathcal{U})$, one has $u\in\Psi(v)$.

(2) If $\mathcal{V}$ is a nontrivial linearity set of $\varTheta_P$, then $\mathcal{U}$ is convex, and for all $v\in \operatorname{cl}(\mathcal{V}) $ and $u\in \Psi(\mathcal{V})$, one has $v\in\Phi(u)$. \upshape

\prof Let us prove the second claim. By Definition \ref{djkyi1}, for any $v^1,v^2\in \mathcal{V}$, one has $\Psi(v^1)=\Psi(v^2)$. From Corollary \ref{interior2}, for any $u\in \Psi(v^1)$, one has $v^1,v^2\in \Phi(u)$. By Corollary \ref{xinzd3}, for any $\alpha\in [0,1]$, one has $v^{\alpha}=\alpha v^1+(1-\alpha)v^2\in \Phi(u)$ and $u\in \Psi(v^{\alpha})$. That is, $\Psi(v^{\alpha})=\Psi(v^1)$. Then $\mathcal{V}$ is convex. The rest of the second claim follows from the above proof trivially. $\square$

In the LP, the actual invariancy region is convex (see Ghaffari-Hadigheh \cite{GG08}). However, in SDP, the actual invariancy region could not be convex. For example, there are six different invariancy sets of $\varTheta_D$ in Example \ref{patex2}. Among these, two open nonlinearity sets are not convex, although the other four linearity sets $\varTheta_D^1,\varTheta_D^2,\varTheta_D^3$ and $\varTheta_D^4$ are convex.

 \cor\label{tplog3} Suppose that $(d,c)\in \operatorname{int}(K \times K^*)$.

(1) Let $\mathcal{U}_1$ and $\mathcal{U}_2$ be two different nontrivial linearity sets of $\varTheta_D$. If $\mathcal{U}=\operatorname{cl}(\mathcal{U}_1)\cap \operatorname{cl}(\mathcal{U}_2)\ne\emptyset$, then $\mathcal{U}$ is a transition face of $\varTheta_D$ and
\[\Phi(\mathcal{U})= \operatorname{conv}(\Phi(\mathcal{U}_1)\cup \Phi(\mathcal{U}_2)). \]

(2) Let $\mathcal{V}_1$ and $\mathcal{V}_2$ be two different nontrivial linearity sets of $\varTheta_P$. If $\mathcal{V}=\operatorname{cl}(\mathcal{V}_1)\cap \operatorname{cl}(\mathcal{V}_2)\ne\emptyset$, then $\mathcal{V}$ is a transition face of $\varTheta_P$ and
\[\Phi(\mathcal{V})= \operatorname{conv}(\Psi(\mathcal{V}_1)\cup \Psi(\mathcal{V}_2)). \] \upshape

\prof By Theorem \ref{zyjg1}, the affine dimensional of the set $\mathcal{U}$ is less than $r$. From Corollary \ref{tplog1}, for any $v^1\in \Phi( \mathcal{U}_1)$, $v^2\in \Phi(\mathcal{U}_2)$, and $u\in\mathcal{U}$, we have
 \[ u \in \Psi(v^1)\cap \Psi(v^2) \quad \operatorname{and} \quad v^1,v^2\in \Phi(u). \] Then by Corollary \ref{xinzd3}, for any $\alpha\in [0,1]$, we have $v^{\alpha}=\alpha v^1+(1-\alpha)v^2\in \Phi(u)$ and $u\in \Psi(v^{\alpha})$. Thus, the first claim is proved. $\square$

\cor\label{tplog4} Suppose that $(d,c)\in \operatorname{int}(K \times K^*)$.

(1) Let $\mathcal{U}_1$ and $\mathcal{U}_2$ be two different nontrivial linearity sets of $\varTheta_D$. If $\mathcal{V}=\Phi(\mathcal{U}_1)\cap \Phi(\mathcal{U}_2)\ne\emptyset$, then $\mathcal{V}$ is a transition face of $\varTheta_P$ and
\[\Psi(\mathcal{V})=\operatorname{cl}(\operatorname{conv}(\mathcal{U}_1\cup \mathcal{U}_2)).\]

(2) Let $\mathcal{V}_1$ and $\mathcal{V}_2$ be two different nontrivial linearity sets of $\varTheta_P$. If $\mathcal{U}=\Psi(\mathcal{V}_1)\cap \Psi(\mathcal{V}_2)\ne\emptyset$, then $\mathcal{U}$ is a transition face of $\varTheta_D$ and
\[\Phi(\mathcal{U})=\operatorname{cl}(\operatorname{conv}(\mathcal{V}_1\cup \mathcal{V}_2)).\]
 \cor\label{tplo92} Suppose that $(d,c)\in \operatorname{int}(K \times K^*)$.

(1) Let $\mathcal{U}$ be a nontrivial linearity set of $\varTheta_D$. If $\mathcal{U}$ contains a recession direction, then $\Phi(\mathcal{U})$ is a transition face of $\varTheta_P$.

(2) Let $\mathcal{V}$ be a nontrivial linearity set of $\varTheta_P$. If $\mathcal{V}$ contains a recession direction, then $\Psi(\mathcal{V})$ is a transition face of $\varTheta_D$. \upshape

The proofs of these two results are similar and have been omitted.

As in Example \ref{xxgh1}, if two vertices $x^*(u^1)$ and $x^*(u^2)$ are adjacent, then $\Phi(u^1)\cap \Phi(u^2)$ is a transition point of $\varTheta_P$. The same result holds for $\varTheta_D$ as well. In Example \ref{patex2}, if $v^1=(1,-1)^T$ and $v^2=(-1,1)^T$, then $\Psi(v^1)\cap\Psi(v^2)=\varTheta_D^3\cap\varTheta_D^4=\{(0,0)^T\}$ is a singleton set. It is easy to verify that $\Phi((0,0)^T)=\operatorname{cl}(\varTheta_P^l)=\operatorname{cl(conv(}\{v^1,v^2\})$ and its origin is the only transition point of $\varTheta_P$.

\subsection{On the existence of a nonlinearity set}
In LP, every slackness vector corresponding to a transition point is a vertex of the polyhedron (for example, see Example \ref{xxgh1}). However, in SDP, the slackness matrix corresponding to a transition point can not be a vertex of the spectrahedon. For example, in Example \ref{patex2}, the origin is only one transition point of $\varTheta_P$, but it is not a vertex. However, the slackness matrix of a transition point is closely related to the vertices. Because the number of vertices of the polyhedron is finite, we make the following conjecture:

{\bf Conjecture:} Number of the vertices of $\mathscr{X}$ or number of transition points of $\varTheta_P$.

An example by Laurent and Poljak \cite{LP95,LP96} confirms that the conjecture is reasonable.

Let $x^*(u^1),x^*(u^2), \cdots,x^*(u^k)$ denote all vertices of $\mathscr{X}$ corresponding to the transition points $u^1,u^2,\cdots,u^k\in\varTheta_D$. We now assume that linear segments
\[ [x^*(u^1),x^*(u^k)], [x^*(u^2),x^*(u^k)],\cdots,[x^*(u^{k-1}),x^*(u^k)] \] do not lie within the boundary of the feasible set of (\ref{primaljihe1}), that is, $\Phi(u^k) \cap \Phi(u^i)=\emptyset$, $i=1,2,\cdots,k-1$. According to Corollaries \ref{mainthe1} and \ref{xinzd3}, the following set
\[ \varTheta_P - \bigcup\limits_{i=1}^k\Phi(u^i)\] is a nonempty open set. Therefore, according to Theorem \ref{zyjg1}, there is a nonlinearity region of $\varTheta_P$. In other words, if there is no linear segment connecting $x^*(u^k)$ and $x^*(u^j) (j=1,2,\cdots,k-1)$ on the boundary of $\mathscr{X}$, then a nonlinearity region exists. Of course, if there is no vertex on the boundary of $\mathscr{X}$, then a nonlinearity region exists. This discussion yields the following results.

\thm Suppose that $(d,c)\in \operatorname{int}(K \times K^*)$. A nonlinearity region exists if one of the following holds:

(1) There is no any vertex on the boundary of the feasible set.

(2) None of the supporting hyperplanes of the feasible set passing through the point $x^*(\bar{u})$ contains other vertex expect for $x^*(\bar{u})$, where $\bar{u}$ denotes a transition point. \upshape

\subsection{The first analysis of multiparametric objective functions}
In this subsection, we discuss the behavior of the multiparametric objective function values in (\ref{primaljihe1}) and (\ref{primaljihe2}).

Consider $u\in \operatorname{int}(\varTheta_D)$, $h\in\mathbb{R}^r$, and $x^*(u+th)$ as $t\rightarrow 0+$. Because the points $x^*(u+th)$ for $t\in(0,\delta)$, where $\delta$ is a positive number, lies in a compact set, $x^*(u+th)$ has a limit point as $t\rightarrow 0+$. Let $x^*_h(u)$ be the limit point. It should be noted that that because of the possible multiplicity of the solutions in $\mathscr{X}^*(u+th)$ and $\mathscr{X}^*(u)$, we shall assume that $x^*(u+th)\rightarrow x^*_h(u)$, because if this is not the case, we choose an appropriate sequence that converges. It is clear that the $x^*_h(u)\in \mathscr{X}^*(u)$.

Analogously, we may define the limit point $y^*_h(v)$ of $y^*(v+th)$ as $t\rightarrow 0+$ and assume that $y^*(v+th)\rightarrow y^*_h(v)$ as $t\rightarrow 0+$.

\lem\label{Gorb1} Let $h\in \mathbb{R}^r$ be arbitrary. Then for any $u\in \operatorname{int}(\varTheta_D)$
\[\lim\limits_{t\rightarrow 0+}\frac{1}{t}\langle c+M^Tu,x^*(u+th)-x^*_h(u)\rangle=0\] and for any $v\in \operatorname{int}(\varTheta_P))$ \[\lim\limits_{t\rightarrow 0+}\frac{1}{t}\langle d+M^Tv, y^*(u+th)-y^*_h(v)\rangle=0. \] \upshape

\prof We follow the proof of Lemma 3.1 in \cite{GS99}. Let us assume that \[ \liminf\limits_{t\rightarrow 0+}\frac{1}{t}\langle c+M^Tu,x^*(u+th)-x^*_h(u)\rangle \leq \varepsilon<0\] (including $\liminf\limits_{t\rightarrow 0+}(\cdot)=-\infty$). Then, there exists a sequence $t_k\rightarrow 0+$ such that
\begin{eqnarray*} && \langle c+M^Tu, x^*(u+t_kh) \rangle\\ & \leq & \langle c+M^Tu, x^*_h(u) \rangle + \varepsilon t_k +o(t_k)\\ &<& \langle c+M^Tu, x^*_h(u) \rangle\end{eqnarray*} for $t_k$ sufficiently small, which contradicts the fact that $x^*_h(u)\in \mathscr{X}^*(u)$.

Similarly, assume that \[ \limsup\limits_{t\rightarrow 0+}\frac{1}{t} \langle c+M^Tu,x^*(u+th)-x^*_h(u)\rangle \geq \varepsilon>0\] (including the case $\liminf\limits_{t\rightarrow 0+}(\cdot)=+\infty$). Then, there exists a sequence $t_k\rightarrow 0+$ such that
\begin{eqnarray*} && \langle d+M^T(u+t_kh), x^*(u+t_kh) \rangle\\ & \geq & \langle d+M^T(u+t_kh), x^*_h(u) \rangle + \varepsilon t_k +t_k\langle M^Tv, x^*(u+th)-x^*_h(u)\rangle +o(t_k). \end{eqnarray*} Since $x^*(u+th)-x^*_h(u)\rightarrow 0$ as $t_k\rightarrow 0+$, it follows from that for $t_k$ sufficiently small, \[ \langle d+M^T(u+t_kh), x^*(u+t_kh) \rangle > \langle d+M^T(u+t_kh), x^*_h(u) \rangle, \] this contradicts the fact that $x^*_h(u)\in \mathscr{X}^*(u+t_kh)$.

Consequently,
\[\liminf\limits_{t\rightarrow 0+}\frac{1}{t} \langle c+M^Tu,x^*(u+th)-x^*_h(u)\rangle=\limsup\limits_{t\rightarrow 0+}\frac{1}{t} \langle c+M^Tu,x^*(u+th)-x^*_h(u)\rangle=0\] and it follows that holds in the statement of the lemma. The second limit is proven in an analogous fashion. $\square$

\thm\label{Gorbthm1} The solution of the following problem produces the directional derivative of $p^*(\cdot)$ at $u\in \operatorname{int}(\varTheta_D)$ in a direction $h\in\mathbb{R}^r$
\begin{equation} \label{gorbg1} p^{*'}(u,h) = \min\limits_{v}\{\langle h,Md+v\rangle| v\in\Phi(u)\}. \end{equation} The solution of the following problem produces the directional derivative of $d^*(\cdot)$ at $v\in \operatorname{int}(\varTheta_P)$ in a direction $h\in\mathbb{R}^r$
\begin{equation} \label{gorbg2} d^{*'}(v,h) = \min\limits_{u}\{\langle h,Mc+u\rangle| u\in\Psi(v)\}. \end{equation} \upshape

\prof Let us now consider the directional derivative of the objective value function $p^*(\cdot)$ at $u\in \operatorname{int}(\varTheta_D)$ in a direction $h\in\mathbb{R}^r$:
\begin{eqnarray*} && \frac{1}{t}(p^*(u+th)-p^*(u))\\ &=&\frac{1}{t} (\langle c+M^T(u+th),x^*(u+th)\rangle -\langle c+M^Tu,x^*(u)\rangle) \\ &=& \langle M^Th,x^*(u+th)\rangle+\frac{1}{t} \langle c+M^Tu,x^*(u+th)-x^*(u)\rangle. \end{eqnarray*}
Then, from Lemma \ref{Gorb1} and $\langle c+M^Tu,x^*(u)\rangle=\langle c+M^Tu,x^*_h(u)\rangle$, we have
\[ p^{*'}(u,h)= \langle M^Th,x^*_h(u)\rangle. \] Because $u\in \operatorname{int}(\varTheta_D)$, $\Phi(u)$ is well defined, that is, $\exists v\in\Phi(u)$, such that
\[\langle M^Th, x^*_h(u)\rangle= \langle h,Mx^*_h(u)\rangle = \langle h,Md+v\rangle \] such that \begin{eqnarray*} p^{*'}(u,h) &=& \min\limits_{x}\{\langle M^Th,x\rangle| Ax=Ad, Mx=Md+v, v\in\Phi(u),x\in K\}\\ & =& \min\limits_{v}\{\langle h,Md+v\rangle| v\in\Phi(u)\}. \end{eqnarray*} The proof is then completed. $\square$

In \cite[Formula 5.227]{BS00}, Bonnans and Shapiro presented a formula for evaluating the derivative of the optimal objection function with respect to the vectors of parameters by calculating the partial derivative of the Lagrangian function. The simpler formulas (\ref{gorbg1}) and (\ref{gorbg2}) can also be viewed as a direct consequence of this work, and cover a similar result for SDP with a scalar parameter (Ref. \cite{GS99}).

\exam {\upshape (Example \ref{patex2} continued).} (1) For the transition point $v^1=(1,1)^T$ of $\varTheta_P$, we have
\[ \left\langle h,\left(\begin{array}{c} c\bullet m^1\\ c\bullet m^2\end{array}\right)+u \right\rangle =u_1h_1+u_2h_2, \] in which $(u_1,u_2)^T\in \Psi(v^1)=\varTheta_D^1$. Then
\[d^{*'}(v^1,h)=\left\{\begin{array}{ll}-h_1-h_2, & \operatorname{if}\ h_1\leq 0\ \operatorname{and} \ h_2\leq 0, \\ -\infty, & \operatorname{if} \ h_1>0\ \operatorname{or}\ h_2>0. \end{array} \right.\]
For the other three transition points of $\varTheta_P$, similar results are also obtained.

For every $v\in\operatorname{int}(\varTheta_P)$, the G$\hat{a}$teaux derivative of $d^*(\cdot)$ at $v$ is equal to
\[d^{*'}(v)=\Psi(v).\]

(2) For the transition point $u^0=(0,0)^T$ of $\varTheta_D$, we have
\[ \left\langle h,\left(\begin{array}{c} d\bullet m^1\\ d\bullet m^2\end{array}\right)+v \right\rangle =v_1h_1+v_2h_2, \] in which $(v_1,v_2)^T\in \Phi(u^0)=\{(v_1,-v_2)^T| -1\leq v_1\ne v_2 \leq 1\}$. Then
\[p^{*'}(u^0,h)=\left\{\begin{array}{ll}h_1-h_2, & \operatorname{if} \ h_1\geq 0\ \operatorname{and} \ h_2\leq 0, \\ -h_1+h_2, & \operatorname{if}\ h_1\leq 0\ \operatorname{and} \ h_2\geq 0, \\ -\infty, & \operatorname{otherwise}. \end{array} \right.\]
For every $u\in \varTheta_D^0$, the G$\hat{a}$teaux derivative of $p^*(\cdot)$ at $u$ is equal to
\[p^{*'}(u)=\Phi(u).\]

\section{Conclusions}
In this study, we propose a lift-and-project procedure to solve a pair of almost primal and dual mpCLOs. We believe that this is one of the first works that focuses
on developing effective tools for the multiparametric analysis of CLOs. Rather than using the optimal basic partition technique for LP or the rank comparison technique for SDP, we show how to use the set-valued mappings to obtain the optimal partition to arbitrary conic representable sets. Similar to the special cases for LP and SDP, it is possible to perform a better parametric analysis based on the optimal partition for perturbations of both the right-hand side and objective function simultaneously.

This paper presented some existing results for a transition point and a nonlinear region, partially answering the open question proposed by Hauenstein et al. \cite{Hau19}. These results depend entirely on the conjecture presented in Section 5. This conjecture is also helpful for understanding the geometry of a conic representable set (Refs \cite{DT15,LP95,LP96}).

Several examples in this study demonstrate the success of the lift-and-project procedure. Currently, we are investigating additional potential applications of CLO using this procedure.

\appendix
\section{Calculations of Example \ref{patex2}}
In this appendix, we evaluate the six indeterminate entries $x_{11},x_{12},x_{23}$ and $\bar{y}_{11}, \bar{y}_{22}, \bar{y}_{33}$ in Example \ref{patex2}. We consider the following two cases:

Case I. The rank of $x^*(u)$ is equal to 1. This condition implies that indeterminate entry $x_{11},x_{12},x_{23}$ must satisfy
\[ \frac{1}{x_{12}} = \frac{x_{12}}{1},\quad \frac{1}{x_{13}} = \frac{x_{13}}{1},\quad \frac{1}{x_{23}} = \frac{x_{23}}{1}, \] which results in $x^*(u)$ being one of the four matrices
\[ \left(\begin{array} {ccc} 1 \\ 1 \\ -1 \end{array} \right)\left(\begin{array} {ccc} 1 \\ 1 \\ -1 \end{array} \right)^T,
\left(\begin{array} {ccc} 1 \\ -1 \\ 1 \end{array} \right)\left(\begin{array} {ccc} 1 \\ -1 \\ 1 \end{array} \right)^T, \left(\begin{array} {ccc} 1 \\ 1 \\ 1 \end{array} \right)\left(\begin{array} {ccc} 1 \\ 1 \\ 1 \end{array} \right)^T, \left(\begin{array} {ccc} 1 \\ -1 \\ -1 \end{array} \right)\left(\begin{array} {ccc} 1 \\ -1 \\ -1 \end{array} \right)^T. \] These matrices are vertices on the surface of the primal feasible set, as shown in Figure \ref{fig:2} (see also the paper \cite{LP96}). By applying the complementary slackness property (\ref{slack1}) to the first vertex, we have
\[\left(\begin{array} {ccc} \bar{y}_{11} & u_1 & -u_2 \\ u_1 & \bar{y}_{22} & -1 \\ -u_2 &-1 & \bar{y}_{33} \end{array} \right) \left(\begin{array} {ccc} 1 \\ 1 \\ -1 \end{array} \right)=\left(\begin{array} {ccc} 0 \\ 0 \\ 0 \end{array} \right). \] By solving this system of linear equations, $\bar{y}_{11},\bar{y}_{22},\bar{y}_{33}$ can be expressed as a function of $(u_1,u_2)\in\mathbb{R}^2$. That is,
\[x^*(u)=\left(\begin{array} {ccc} 1&1&-1 \\ 1&1&-1 \\ -1 &-1&1 \end{array} \right),\quad \bar{y}^*(u)=\left(\begin{array} {ccc} -u_1-u_2 & u_1 & -u_2 \\ u_1 & -u_1-1 & -1 \\ -u_2 &-1 & -u_2-1 \end{array} \right). \] As $\bar{y}^*(u)$ is positive semidefinite, it is easy to verify that $u=(u_1,u_2)^T$ belongs to the following set
 \[\varTheta_D^1 =\{(u_1,u_2)^T |u_1< -1, u_2< -1, u_1+u_2+u_1u_2\geq0 \}. \]
This region is a convex set, bounded by one branch of a hyperbola in plane. Hence for any $u=(u_1,u_2)^T\in\varTheta_D^1$, one has \[\Phi(u)=\frac{1}{2}\left(\begin{array} {ccc} m^1\bullet x^*(u) \\ m^2\bullet x^*(u) \end{array} \right)=
\left(\begin{array} {ccc} 1 \\ 1 \end{array} \right) \] since $ m^i\bullet (x^*(u)-d)=2v_i$ for $i=1,2$. Moreover, one has $\Psi((1,1)^T)=\varTheta_D^1$.

Analogously, we also obtain all linearity sets associated with the other three vertices. Detailed results are presented in Table 2.

Case II: The rank of $\bar{y}^*(u)$ is equal to 1. It follows from
\[ \frac{\bar{y}_{11}}{u_1} = \frac{-u_2}{-1}, \quad \frac{u_1}{-u_2}= \frac{\bar{y}_{22}}{-1}=\frac{-1}{\bar{y}_{33}} \] that $\bar{y}^*(u)$ can be expressed in terms of the parameters $u_1$ and $u_2$ as follows:
\[\bar{y}^*(u)=\left(\begin{array} {ccc} u_1u_2 & u_1 & -u_2 \\ u_1 & \frac{u_1}{u_2} & -1 \\ -u_2 &-1 & \frac{u_2}{u_1} \end{array} \right)= \frac{1}{u_1u_2} \left(\begin{array} {ccc} u_1u_2 \\ u_1 \\ -u_2 \end{array} \right) \left(\begin{array} {ccc} u_1u_2 \\ u_1 \\ -u_2 \end{array} \right)^T,\]
where $u_1u_2>0$. It follows from
\[ \left(\begin{array} {ccc} 1 & x_{12} & x_{13} \\ x_{12} & 1 & x_{23} \\ x_{13}& x_{23} & 1\end{array} \right) \left(\begin{array} {ccc} u_1u_2 \\ u_1 \\ -u_2 \end{array} \right)=\left(\begin{array} {ccc} 0 \\ 0 \\ 0 \end{array} \right), \] that
\[ \begin{array} {ccc} x_{12} & = & \frac{u_2}{2u_1^2} - \frac{u_2}{2} -\frac{1}{2u_2}, \\ x_{13} &=& \frac{1}{2u_1} + \frac{u_1}{2} -\frac{u_1}{2u_2^2}, \\ x_{23} &=& \frac{u_2}{2u_1} +\frac{u_1}{2u_2}-\frac{u_1u_2}{2}. \end{array} \] Thus, six indeterminate entries are established. In particular, the case of $u_1=u_2$ was discussed by Mohammad-Nezhad and Terlaky \cite{MT20}.

Note that $x^*(u)$ is positive semidefinite if and only if $|x_{23}|\leq 1$, i.e., \[ -2\leq \frac{u_2}{u_1} +\frac{u_1}{u_2}-u_1u_2\leq 2.\] Then $u_1u_2>0$ yields that
\[-2u_1u_2\leq u_1^2+u_2^2-u_1^2u_2^2\leq 2u_1u_2,\] or equivalently,
\[0\leq (u_1+u_2)^2-(u_1u_2)^2,\qquad and \qquad (u_1-u_2)^2-(u_1u_2)^2\leq 0.\] This concludes the following inequalities
\begin{eqnarray} \label{bdsbuj1} (u_1+u_2+u_1u_2)(u_1+u_2-u_1u_2) &\geq& 0, \\ \label{bdsbuj2} (u_1-u_2+u_1u_2)(u_1-u_2-u_1u_2)&\leq& 0. \end{eqnarray}

\begin{figure} \centering \includegraphics[height=4cm,width=5cm]{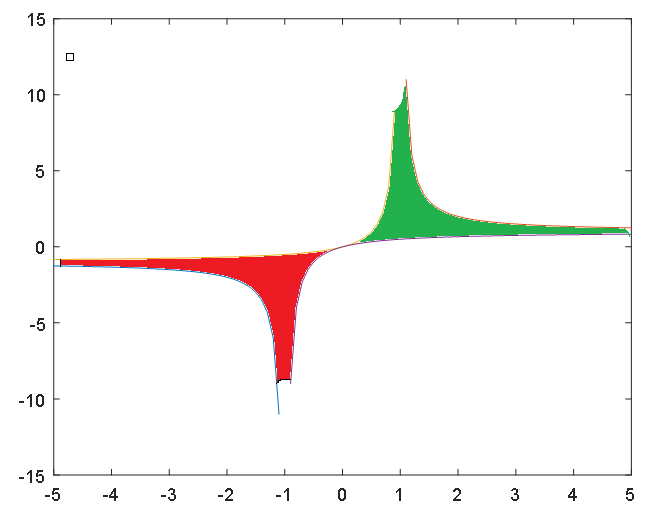} \caption{Primal conic representable set separated by four curves} \label{fig:3} \end{figure}

\lem\label{quyu1} Define a set
\[\varTheta_D^0=\{(u_1,u_2)^T\in \mathbb{R}^2| u_1u_2>0, \operatorname{ strict \ inequalities\ (\ref{bdsbuj1}) \ and \ (\ref{bdsbuj2})\ hold} \}. \] Then $\varTheta_D^0= \mathbb{R}^2 -\bigcup\limits_{i=1}^4\varTheta_D^i$. \upshape

\prof Define four curves as follow
\[ \begin{array}{lll} l_1: & u_1+u_2+u_1u_2=0, & u_1<-1,u_2<-1,\\
l_2: & -u_1-u_2+u_1u_2=0, & u_1>1,u_2>1, \\ l_3: & u_2-u_1-u_1u_2=0, & u_1<1,u_2>-1, \\ l_4: & u_1-u_2-u_1u_2=0, & u_1>-1,u_2<1. \end{array} \] Each of them represents a unilateral branch of the hyperbola, in which $l_1$ and $l_2$ are symmetric about line $u_1+u_2=0$, and $l_3$ and $l_4$ is symmetric about line $u_1-u_2=0$. These four curves define some areas in $\mathbb{R}^2$. Instead of a formal and tedious proof we plot these areas in Figure \ref{fig:3}. The four outer closed regions are $\varTheta_D^1$, $\varTheta_D^2$, $\varTheta_D^3$ and $\varTheta_D^4$. Intuitively, $\varTheta_D^0$ lies in the two curved triangles inside (the colored part in Figure \ref{fig:3}), where the two strict inequalities (\ref{bdsbuj1}) and (\ref{bdsbuj2}) hold for any $(u_1,u_2)\in\{(u_1,u_2)\in\mathbb{R}^2|u_1u_2>0\}$. $\square$

\lem\label{quyu2} For any $u\in \varTheta_D^0$, we have $|x_{12}|< 1$ and $|x_{13}|< 1$. \upshape

\prof Assume that $u_1>0$ and $u_2>0$. Then, $\varTheta_D^0$ in the first quadrant is bounded by three curves $l_2,l_3,l_4$. That is,
\begin{eqnarray*} u_1+u_2-u_1u_2&>&0, \\
u_2-u_1-u_1u_2&<&0, \\ u_1-u_2-u_1u_2&>&0, \\ u_1>0,\ u_2&>&0. \end{eqnarray*} Or equivalently,
\begin{eqnarray*} 0<&\frac{1}{u_1}+\frac{1}{u_2}&<1, \\
-1<&\frac{1}{u_1}-\frac{1}{u_2}&<1. \end{eqnarray*} Therefore, we have \[ \frac{1}{u_1^2}-\frac{1}{u_2^2} <1<1+ \frac{2}{u_2}, \]
that is, $x_{12}< 1$. On the other hand, we have
 \[2x_{12}= \left(\frac{1}{u_1}+\frac{u_2}{u_1^2}-\frac{u_2}{u_1}\right)+
 \left(1-\frac{1}{u_1}-\frac{1}{u_2}\right)+\left(\frac{u_2}{u_1}-u_2+1\right)-2>-2. \] The first claim for $u_1<0$ and $u_2<0$ and the second claim are proven in an analogous fashion. $\square$

From Lemma \ref{quyu1}, for any $u=(u_1,u_2)^T \in \varTheta_D^0$, there is a pair of optimal solution $(x^*(u),\bar{y}^*(u))$ such that the ranks of $x^*(u)$ and $\bar{y}^*(u)$ are equal to two and one, respectively. Moreover, one has \begin{equation} \label{setvalue11} \Phi(u)=\frac{1}{2}\left( \begin{array} {ccc} m^1\bullet (x^*(u)-d) \\ m^2\bullet (x^*(u)-d) \end{array} \right) = \left( \begin{array} {ccc} \frac{u_2}{2u_1^2} - \frac{u_2}{2} -\frac{1}{2u_2} \\ -\frac{1}{2u_1} -\frac{u_1}{2} +\frac{u_1}{2u_2^2} \end{array} \right). \end{equation}

Similarly, let
\[ \bar{x}^*(v)=\left(\begin{array}{ccc} 1&v_1&-v_2\\ v_1 & 1 &\bar{x}_{23} \\ -v_2 & \bar{x}_{23} & 1 \end{array}\right), \quad y^*(v)=\left(\begin{array}{ccc} y_1\\ y_2\\ y_3\end{array}\right)^T\left(\begin{array}{ccc} y_1\\ y_2\\ y_3\end{array}\right) \] denote the optimal solutions of problem (\ref{mostl2}) and its nonstandard dual, in which $det(x^*(v))=0$ and $y_2y_3=-1$ are assumed. The first assumption implies that
\[1 - 2v_1v_2 \bar{x}_{23} - v_1^2-v_2^2-\bar{x}_{23}^2=0. \] Solving the quadratic equation with one variable, $\bar{x}_{23}$, to obtain
\[\bar{x}_{23}=v_1v_2\pm \sqrt{(1-v_1^2)(1-v_2^2)}. \] However, it follows from the complement slackness property (\ref{slack2}) that
\[ \left(\begin{array}{ccc} 1&v_1&-v_2\\ v_1 & 1 &\bar{x}_{23} \\ -v_2 & \bar{x}_{23} & 1 \end{array}\right) \left(\begin{array}{ccc} y_1\\ y_2\\ y_3\end{array}\right) =0,\] to yield
\[ (1-v_1^2)y_2=(v_1v_2+\bar{x}_{23})y_3.\] Then the second assumption implies that \[\bar{x}_{23}=v_1v_2-\sqrt{(1-v_1^2)(1-v_2^2)}. \] Therefore, one has
\[(y_1,y_2,y_3)^T =\left(\frac{v_2\sqrt{1-v_2^2}+v_1\sqrt{1-v_1^2}}{\sqrt[4]{(1-v_1^2)(1-v_2^2)}},-\sqrt[4]{\frac{1-v_2^2}{1-v_1^2}}, \sqrt[4]{\frac{1-v_1^2}{1-v_2^2}}\right)^T\] such that
\begin{equation} \label{setvalue12} \Psi(v)= \frac{1}{2} \left(\begin{array}{c} m^1\bullet (y^*(v)-c)\\ m^2\bullet (y^*(v)-c) \end{array} \right)=\left(\begin{array}{c}-v_2-v_1\sqrt{\frac{1-v_2^2}{1-v_1^2}} \\ -v_1-v_2 \sqrt{\frac{1-v_1^2}{1-v_2^2}} \end{array} \right). \end{equation}


\end{document}